\documentclass[cmp]{svjour}
\usepackage{amsmath,amsfonts,amssymb}
\usepackage{graphicx}

\def\G{{\cal G}}
\def\D{{\cal D}}
\def\P{{\cal P}}
\def\Q{{\cal Q}}
\def\lo{\longrightarrow}
\def\ov{\overline}

\newcommand{\spn}{\operatorname{span}}

\newcommand{\Legendre}{{\cal L}}
\newcommand{\until}[1]{\{1,\dots, #1\}}

\newcommand{\real}{{\mathbb{R}}}

\newcommand{\sphere}{{\mathbb{S}}}

\newcommand{\graph}[1]{\text{Graph}(#1)}
\newcommand{\ann}[1]{\text{Ann}(#1)}
\newcommand{\pder}[2]{\frac{\partial #1}{\partial #2}}

\newcommand{\setdef}[2]{\{#1 \; | \; #2 \}}

\renewcommand{\theenumi}{(\roman{enumi})}

\renewcommand{\square}{\Delta}

\begin{document}
\title{Hamiltonian theory of constrained impulsive motion}

\author{Jorge Cort\'es\inst{1} \and Alexandre M. Vinogradov\inst{2}}

\institute{Coordinated Science Laboratory, University of Illinois at
  Urbana-Champaign, 1308 W Main St, Urbana, IL 61801, USA,
  \email{jcortes@uiuc.edu}, Phone: +1-217-244-8734, Fax:
  +1-217-244-1653 \and Dipartimento di Matematica e Informatica,
  Facolt\`a di Scienze Matematiche Fisiche e Naturali, Universit\`a di
  Salerno, Via S.  Allende, I-84081 Baronissi, and Istituto Nazionale
  di Fisica Nucleare, sezione Napoli-Salerno, Italy,
  \email{vinograd@unisa.it}}


\date{\today}

\maketitle

\begin{abstract}
  This paper considers systems subject to nonholonomic constraints which
  are not uniform on the whole configuration manifold.  When the
  constraints change, the system undergoes a transition in order to comply
  with the new imposed conditions.  Building on previous work on the
  Hamiltonian theory of impact, we tackle the problem of mathematically
  describing the classes of transitions that can occur. We propose a
  comprehensive formulation of the Transition Principle that encompasses
  the various impulsive regimes of Hamiltonian systems.  Our formulation is
  based on the partial symplectic formalism, which provides a suitable
  framework for the dynamics of nonholonomic systems.  We pay special
  attention to mechanical systems and illustrate the results with several
  examples.
\end{abstract}



\section{Introduction}

In this paper we consider the problem of mathematically describing
impulsive motions (impacts, collisions, reflection, refractions) of
Hamiltonian systems subject to nonholonomic constraints. An impulsive
behavior takes place when one or more of the basic ingredients of the
Hamiltonian dynamical picture undergoes a drastic change. As an example,
one may consider the instant of time when the configuration space of the
system collapses instantaneously because of an inelastic collision. Another
example is given by a ray of light that splits into reflected and refracted
rays when passing from one optical media to another, and so on. In
situations like these, the phase trajectory of the system becomes
discontinuous and the problem of how to describe this discontinuity arises.

The problem of describing impulsive motion has been extensively studied in
classical books such as~\cite{Ap,Mo,Pa,Par,Ro}. In these references, the
emphasis is put on the analysis of mechanical systems subject to impulsive
forces, and in particular, the study of rigid body collisions by means of
Newton and Poisson laws of impact.  Impulsive nonholonomic constraints
(i.e. constraints whose reaction force is impulsive) are also considered
in~\cite{NeFu,Ro}, and, from a geometric perspective, in more recent works
such as~\cite{CeLaRe,IbLeLaMaMaPi,LaTu}. If impulsive constraints and
impulsive forces are present at the same time, Newton and Poisson
approaches have been revealed to be physically inconsistent in certain
cases~\cite{Br,St}.  This surprising consequence of the impact laws is only
present when the velocity along the impact surface stops or reverses during
collision due to the friction.  Energetically consistent hypothesis for
rigid body collisions with slip and friction are proposed in~\cite{Ste,St}.
From a design point of view, the interest in systems subject to impulsive
forces is linked to the emergence of nonsmooth and hybrid dynamical systems
in Control Theory, i.e., systems where continuous and discrete dynamics
coexist, see~\cite{Bro,Br,Br2,ScSc} and references therein.  Hybrid
mechanical systems that locomote by switching between different constraint
regimes and are subject to elastic impacts are studied in~\cite{BuZe}.
Hyper-impulsive control of mechanical systems is analyzed in~\cite{Fa}.

Here, we aim at a comprehensive analysis of the various situations
which can occur concerning impulsive regimes of nonholonomic
Hamiltonian systems.  In particular, we focus on two different but
complementary cases.  The first one deals with a drastic change in the
nonholonomic constraints imposed on the system.  The second one
concerns a drastic change of the Hamiltonian function and includes, in
particular, collisions and impacts of nonholonomic systems. The
proposed solution is given in terms of a generalized version of the
\emph{Transition Principle}.  This principle, sketched for the first
time in a series of lectures of the second author~\cite{BoVi} for
discontinuous Hamiltonians, was recently extended to other
non-constrained situations in~\cite{PuVi,PuVi2} (see also~\cite{CoMa}
for a related discussion in an optimal control setting). By its very
nature, the Transition Principle is a direct dynamic interpretation of
the geometric data of the problem.  This feature distinguishes it from
other approaches.  For instance, in Classical Mechanics, the velocity
jumps caused by an impact are traditionally derived from some
assumptions on the nature of the impulsive forces (see, for
instance,~\cite{Ap,LeAm}).  However, these assumptions are not logical
consequences of the fundamental dynamical principles and therefore one
should really consider them as additional principles for impulsive
problems.
The distinguishing feature of the Transition Principle is that it
gives full credit to the geometry of the nonholonomic Hamiltonian
system.  This seems reasonable for the impulsive regime, keeping in
mind the perfect adequacy of the Hamiltonian description to the
dynamical behavior of the system in the absence of impulsive motions.
In addition, there are some noticeable advantages deriving from the
Transition Principle.  First of all, its application gives an exact
and direct description of the post-impact state which is of immediate
use for both theoretical and computational purposes. Secondly, it is
still valid in some situations where standard methods can be hardly
applied.  In particular, this is the case of Hamiltonian systems
describing the propagation of singularities of solutions of partial
differential equations (consider, for instance, the example of
geometrical optics)~\cite{LiMaSpVi,Ly,Vi}.  Clearly, no variational or
traditional approach can be applied to this very important class of
systems.

A second contribution of this paper concerns the formulation of the
dynamics of nonholonomic Hamiltonian systems. We make use of the notion of
partial symplectic structures introduced in~\cite{BoVi} and relate this
framework with other modern approaches to nonholonomic systems
(see~\cite{BaSn,Bl,BlKrMaMu,Co,Ko,LeMa,Mar,ScMa} and references therein).
One advantage of the partial symplectic formalism is that it allows us to
draw clear analogies between the unconstrained and constrained situations.
Another advantage is that the treatment of nonlinear constraints can be
easily incorporated.

The paper is organized as follows. Section~\ref{se:preliminaries}
introduces some geometric preliminaries on distributions, constraint
submanifolds and partial symplectic structures.  In
Section~\ref{se:dynamics-nonholonomic}, we show how any nonholonomic
Hamiltonian system possesses an associated partial symplectic
structure, and we use this fact to intrinsically formulate the
dynamics. We also analyze systems with instantaneous nonholonomic
constraints and systems exhibiting discontinuities.  In
Section~\ref{se:Transition-Principle}, we develop a new formulation of
the Transition Principle for systems with constraints.  We present the
novel notion of focusing points with respect to a constraint
submanifold and we also introduce the concept of constrained
characteristic. Decisive points are defined for each impulsive regime
resorting to in, out and trapping points.
Section~\ref{se:mechanical-systems} presents a detailed study of the
concepts introduced in the previous sections in the case of mechanical
systems. We compute the focusing points and the characteristic curves,
and present various results concerning the decisive points. We also
prove an appropriate version for generalized constraints of the
classical Carnot's theorem for systems subject to impulsive forces: if
the constraints are linear, we show that the Transition Principle
always implies a loss of energy. We conclude this section by showing
that if the constraints are integrable, then our formulation of the
Transition Principle recovers the solution for completely inelastic
collisions~\cite{PuVi2}. Finally, Section~\ref{se:examples} presents
various examples of the application of the above-developed theory.

To ease the exposition, below we make use of the standard notation
concerning differential geometry and the Hamiltonian formalism without
making explicit reference to any work. In particular, we denote by
$\Lambda^i$ (resp., $D(M)$) the $C^{\infty}(M)$-module of $i$th order
differential forms (resp., of vector fields) on a manifold $M$.  We
use $F^*(\varphi)$ to denote the pullback with respect to a smooth map
$F$ of a function or differential form $\varphi$.  If $x$ is a point
of $M$, then the subscript $x$ refers to the value of the
corresponding geometric object at $x$. For instance, $X_x$ stands for
the vector field vector $X\in D(M)$ evaluated at $x$.  The interested
reader may consult classical books such as~\cite{AbMaRa,KoNo,LiMa} for
further reference. We also assume smoothness of all the objects we are
dealing with.

\section{Preliminaries}\label{se:preliminaries}

In this paper we deal with Hamiltonian systems defined on the
cotangent bundle $T^*M$ of an $n$-dimensional manifold $M$. In the
particular case of a mechanical system, $M$ and $T^*M$ are,
respectively, the configuration space and the phase space of the
system.  As usual, $\pi_M:T^*M \rightarrow M$ (or simply $\pi$) stands
for the canonical projection from $T^*M$ to $M$, $H\in
C^{\infty}(T^*M)$ for the Hamiltonian function and $X_H \in D(M)$ for
the corresponding Hamiltonian vector field. The canonical symplectic
structure on $T^*M$ is denoted by $\Omega =\Omega_M$. In canonical
coordinates $(q^a,p_a)$, $a = 1, \dots,n$ of $T^*M$, the symplectic
form reads $\Omega=dq^a \wedge dp_a$.

We say that the Hamiltonian system $(M,H)$ \emph{comes from a
  Lagrangian system $(M,L)$ on $TM$} if $H =
(\Legendre^*_L)^{-1}(E_L)$, where $E_L\in C^ {\infty}(TM)$ is the
energy function corresponding to the (hyper-regular) Lagrangian
function $L\in C^{\infty}(TM)$ and $\Legendre_L:TM \rightarrow T^*M$
is the associated Legendre map.

If $X$ is a vector field on $T^*M$, then the map $\alpha_X : T^*M
\rightarrow TM$ defined by
\[
\alpha_X(\theta)=d_{\theta}\pi (X_{\theta})\in T_{\pi(\theta)} M
\, , \quad \theta\in T^*M \, ,
\]
denotes the \emph{anti-Legendre map associated with $X$}. In standard
coordinates, if $X=A^a(q,p)\pder{}{q^a}+B^a(q,p)\pder{}{p_a}$, then
$\alpha_X$ reads $\alpha_X(q^a, p_a)=(q^a, A^a(q,p))$. For the
Hamiltonian vector field $X=X_H$, we write $\alpha_H$ instead of
$\alpha_{X_H}$, so that
\[
\alpha_H : (q,p) \mapsto \left( q, v=\pder{H}{p} \right) \, .
\]
It is not difficult to see that if the Hamiltonian system $(M,H)$
comes from a Lagrangian system $(M,L)$, then
$\alpha_H=\Legendre_L^{-1}$.

\subsection{Distributions and
  codistributions}\label{se:distributions}

Recall that a \emph{distribution} (resp., \emph{codistribution}) on a
manifold $M$ is a vector subbundle of $TM$ (resp., of $T^*M$).  The
\emph{annihilator} of a distribution $\D$ on $M$ is the codistribution
$\ann{\D}$ defined by
\[
\ann{\D}_x = \setdef{\theta \in T_x^*M}{\theta(\xi)=0, \forall \xi\in
  \D_x} \, , \quad x \in M \, .
\]
If $\D$ is $(n-m)$-dimensional, the codistribution $\ann{\D}$ is
$m$-dimensional. The \emph{dual bundle $\D^*$ of $\D$} is canonically
identified with the cotangent bundle $T^*M$ modulo $Ann(\D)$. We will
also denote by $\D^\perp$ the orthogonal complement of a distribution
$\D$ on $T^*M$ with respect to the symplectic form $\Omega$, i.e.,
\[
\D_y^\perp = \setdef{\xi \in T_y(T^*M)}{\Omega_y (\xi,\eta) = 0 \, , \;
  \forall \eta \in \D_y} \, , \quad y\in T^*M \, .
\]

A vector field $X\in D(M)$ \emph{belongs} to $\D$ if $X_x
\in \D_x$ for all $x\in M$. Vector fields belonging to $\D$ constitute a
$C^{\infty}(M)$-module, denoted by $D_{\D}(M)$, which is a
submodule of $D(M)$.  In the partial symplectic formalism
(see Section~\ref{se:partial} below), they are interpreted as
``constrained'' vector fields. Dually, denote by $\Lambda_{\D}^1(M)$
the $C^{\infty}(M)$-module of sections of the bundle $\D^*$ and by
$\Lambda_{\D}^i(M)$ its $i$th exterior product. These are interpreted
as ``constrained'' differential $i$-forms. We denote the natural
restriction map from $\Lambda^i(M)$ to $\Lambda_{\D}^i(M)$ by
$r_{\D}:\Lambda^i(M) \rightarrow \Lambda_{\D}^i(M)$.

The geometric description of nonholonomic systems in the framework of
the partial symplectic formalism~\cite{BoVi} requires a slight
``affine'' generalization of these standard notions. Namely, an
\emph{affine distribution} on a manifold $M$ is an affine subbundle
$\square$ of $TM$. This means that the fiber $\square_x$ of $\square$
over $x\in M$ is an affine subspace in $T_x M$.  Therefore,
$\square_x$ can be represented in the form $\square_x = v +\Delta^0_x$
with $v \in T_x M$ and $\Delta^0_x$ being the vector subspace of $T_x
M$ canonically associated with $\square_x$. In this representation,
the \emph{displacement vector} $v$ is unique modulo $\Delta^0_x$. The
union $\cup_{x \in M}\Delta^0_x$ constitutes a linear distribution of
the tangent bundle $TM$, denoted by $\Delta^0$, canonically associated
with $\square$.  It is not difficult to see that there always exist a
vector field $Y \in D(M)$ such that $Y_x$ is a displacement vector for
$\square_x$. Such vector fields are called \emph{displacement vector
  fields of $\square$}.
Obviously, displacement vector fields differ by another vector field
belonging to $\Delta$.  In coordinate terms, an $(n-m)$-dimensional
affine codistribution is described by a system of linear equations
$\Phi_i=0$ with respect to the variables $p_a$, i.e., $\Phi_i (q,p) =
\Phi_{ia}(q)p_a + \Phi_{i0}(q)$, $i=1,\dots ,m$.

Similarly, an \emph{affine codistribution} on $M$ is an affine
subbundle $C\subset T^*M$ of the cotangent bundle. As above, one has
$C = \Upsilon + C^0$, where $C^0$ is the unique codistribution on $M$
canonically associated with $C$, and $\Upsilon \in \Lambda^1(M)$ is a
\emph {displacement form}.  Point-wisely this means that $C_x =
\Upsilon_x + C^0_x$, for all $x\in M$.

\subsubsection{Linear constraints}\label{se:lconstraints}

In the case of linear constraints, the analogy between free and
constrained systems is particularly clear.  In fact, it is natural to
interpret an affine distribution (resp., codistribution) on a manifold
$M$ as the ``constrained'' tangent (resp., cotangent) bundle of $M$.
A linearly constrained Hamiltonian system is then a triple $(M,H,C)$,
with $H\in C^{\infty}(T^*M)$ and $C$ an affine codistribution on $M$.
Similarly, a triple $(M,L,\Delta)$, with $L\in C^{\infty}(TM)$ and
$\Delta$ an affine distribution on $M$, is a linearly constrained
Lagrangian system.  The anti-Legendre map allows one to pass from a
constrained Hamiltonian system to the corresponding Lagrangian system
and vice versa.  More precisely, if $(M,H,C)$ is a linearly
constrained Hamiltonian system, the map $\alpha_H$ is linear and
$(M,H)$ comes from a Lagrangian system $(M,L)$, then the corresponding
linearly constrained Lagrangian system is $(M,L,\Delta)$, with
$\Delta=\alpha_H(C)$. To go in the opposite direction, one must use
the Legendre map $\Legendre_L$ instead of $\alpha_H$.

Throughout the paper, we distinguish the class of mechanical systems
subject to linear constraints because of two reasons. First,
classically they have been intensively studied.  Second, one can
extract from them the motivations for the basic constructions which
will be discussed below.

\subsubsection{Nonlinear constraints}\label{se:constraints}

In the Hamiltonian setting, the \emph{nonholonomic constraints} are
given by a submanifold (not necessarily a vector subbundle) $C\subset
T^*M$.
Similarly, \emph{nonholonomic Lagrangian or kinematic constraints} are
given by a submanifold $C'\subset TM$. If the Hamiltonian system
$(M,H)$ comes from a Lagrangian system $(M,L)$, then
$C=\Legendre_L(C')$ if and only if $C'=\alpha_H(C)$.  In mechanics,
these two approaches correspond to two possible descriptions of
nonholonomic constraints: either as limitations imposed on the momenta
or as limitations imposed on the velocities, respectively.  The fact
that $C$ (resp., $C'$) represents limitations imposed only on the
momenta (resp., velocities), but not on the configurations of the
system, implies that the projection $\pi$ must send $C$ (resp., $C'$)
surjectively onto $M$.  However, the assumption of ``infinitesimal
surjectivity'' of $\pi|_C$ is more adequate in this context. This
means that $\pi|_C$ is a \emph{submersion}, i.e., $d_y(\pi|_C) : T_yC
\rightarrow T_{\pi(y)} M$ is surjective for all $y\in C$.  With this
motivation, we adopt the following definition.
\begin{definition}
  A set of \emph{nonholonomic constraints} imposed on a Hamiltonian
  system $(M,H)$ is a submanifold $C\subset T^*M$ such that $\pi|_C$
  is a submersion.  The constrained Hamiltonian system is denoted by
  $(M,H,C)$.
\end{definition}
Since $C$ and $\alpha_H(C)$ are, respectively, interpreted as the
constrained cotangent and tangent bundle of the system $(M,H,C)$, we
will always assume that they have equal dimensions.  It is worth
stressing that the above definition makes also sense for manifolds
with boundary. In such a case, the boundary of $T^*M$ is
$\pi^{-1}(\partial M)$ and the boundary of $C$ is~$\partial C = C\cap
\pi^{-1 }(\partial M)$.

\begin{remark}
  {\rm Similarly, \emph{nonholonomic Lagrangian constraints} are
    represented by submanifolds of $TM$ that project regularly onto
    $M$.}
\end{remark}

In what follows, $\Phi_i(q,p)=0$, $i=1,...,m$, will denote a set of
local equations defining $C$. For a point $x\in M$, we denote by $C_x$
the fiber of $C$ at $x$,
\[
C_x = C \cap T_x^*M = \{ y \in C \mid \pi(y) = x \}
=(\pi_{|C})^{-1}(x) \, .
\]

\subsubsection{Instantaneous nonholonomic
  constraints}\label{se:inst-constraints}

Let $N$ be a hypersurface in $M$. Consider the induced hypersurface
$T^*_N M = \pi^{-1} (N)$ of $T^*M$. Let $C \subset T^*M$ be a set of
nonholonomic constraints on $M$.  Instantaneous constraints may be
thought as limitations on the momenta (resp., velocities) of the
system that are imposed only at the instant when a trajectory passes
through a point of $N$.  Therefore, they are represented by a
submanifold~$C^{\text{inst}}$ of~$C\cap T^*_N M$.  These constraints
are assumed to be additional to the ones already prescribed by $C$.
In order to admit an adequate mechanical interpretation, we also
assume that the projection $\pi$ restricted to $C^{\text{inst}}$ is a
submersion onto $N$.  From the Lagrangian point of view, instantaneous
kinematic constraints are naturally interpreted as a submanifold
${C^{\text{inst}}}'$ of $TN$.  Based on these considerations, we take
the following definition.

\begin{definition}\label{dfn:inst-constraints}
  Let $(M,H,C)$ be a constrained Hamiltonian system and let $N$ be a
  hypersurface of $M$. A set of \emph{instantaneous constraints along $N$}
  imposed on $(M,H,C)$ is a submanifold~$C^{\text{inst}}$ of~$C\cap T^*_N
  M$ such that $\pi$ restricted to~$C^{\text{inst}}$ is a submersion onto
  $N$.
\end{definition}

It is worth stressing that, in some cases, a set of instantaneous
constraints along $N$ additionally verifies the condition
$C^{\text{inst}}\subset \alpha_H^{-1}(TN)$ (here $TN$ is thought to be
naturally embedded into $TM$). In an inelastic scenario, where the
nonholonomic motion in $M$ is forced to take place in $N$ after the impact,
this latter condition formalizes the parity between the Hamiltonian and
Lagrangian approaches: if the Hamiltonian system in question comes from a
Lagrangian one, then $C^{\text{inst}} = \Legendre_L({C^{\text{inst}}}')$,
with ${C^{\text{inst}}}'$ being the instantaneous kinematic constraints.

\subsection{Dynamics of Hamiltonian systems}\label{se:dynamics}

As is well-known, in the absence of constraints, the dynamics of the
Hamiltonian system $(M,H)$ is given by the Hamiltonian vector field $X_H$,
whose coordinate description is
\begin{align*}
  \frac{dq^a}{dt} = \pder{H}{p_a} \, , \quad \frac{dp_a}{dt} = -
  \pder{H}{q^a} \, , \quad a = 1, \dots, n \, .
\end{align*}
In the presence of constraints, the ``free'' Hamiltonian vector field
$X_H$ must be modified along the constraint manifold $C$ in order to
become tangent to $C$. In the traditional approach this goal is
achieved by adding to $X_H$ another vector field along $C$, say, $R$,
interpreted as the \emph{reaction of constraints}. From a purely
geometrical point of view, the choice of a vector field that makes
$X_H$ tangent to $C$ is far from being unique.  Therefore, a new
principle must be invoked to select the one that merits to be called
the ``reaction of constraints''. The history of this problem (see, for
instance,~\cite{NeFu}) shows that its solution is not straightforward.
By applying, for instance, the Lagrange-d'Alembert principle
(see~\cite{Bl,Co,NeFu}), one gets the following equations of motion
\begin{align*}
  \frac{dq^a}{dt} = \pder{H}{p_a} \, , \quad \frac{dp_a}{dt} = -
  \pder{H}{q^a} + \lambda_i \pder{\Phi_i}{p_a} \, , \quad \Phi_i
  (q^a,p_a) = 0 \, ,
\end{align*}
$a = 1, \dots, n$, $i = 1, \dots, m$, where the ``Lagrange
multipliers'' $\lambda_i$'s are to be duly determined. Shortcomings of
such an approach are that it is not manifestly intrinsic and does not
reveal clearly the geometric background of the situation.  This is why
in our further exposition we shall follow a purely geometric approach,
which does not require any discussion of reactions of constraints.  It
is based on the concept of \emph{partial symplectic formalism}, which
also appears to be more concise from an algorithmic point of view.

\subsection{Partial symplectic structures}\label{se:partial}

The following elementary facts from linear algebra will be most
useful. Let $V$ be a vector space, $W\subset V$ a subspace and $b:V
\times V \rightarrow \real$ a bilinear form on $V$. Denote by
$W_b^\perp$ the $b$-orthogonal complement of $W$,
\[
W_b^\perp = \setdef{v\in V}{b(v,w)=0, \forall w\in W} \, .
\]
Note that $W\cap W_b^\perp =0$ if and only if the restriction
$b\mid_W$ of $b$ to $W$ is nondegenerate. The form $b$ is said to be
\emph{nondegenerate} on an affine subspace $U$ of $V$, $U=p_0 + W$,
$p_0 \in V$, if it is nondegenerate on its associated vector space
$W$. In such a case, $U$ can be uniquely represented in the form
$U=p_1+W$ with $p_1\in W_b^\perp$ due to the fact that $U\cap
W_b^\perp =\{p_1\}$. The vector $p_1$ is called the \emph{canonical
  displacement} of $U$ with respect to~$b$. Consider the associated
map
\[
\top_{W,b} : W\lo W^*, \qquad \top_{W,b}(w)=b(w,\cdot) \, , \; w \in
W.
\]
In other words, $ \top_{W,b}(w)(w')=b(w,w')$, for all $w'\in W$.
Obviously, $\top_{W,b}$ is an isomorphism if and only if $b_{|W}$ is
nondegenerate.

If $b$ is skew-symmetric and nondegenerate on $V$, and $W$ is a
subspace of $V$ with codimension one, then the kernel of the
restricted form $b|_W$, $\ker b|_W$ is a 1-dimensional subspace, i.e.,
a \emph{line} in $V$ contained in $W$.  Therefore, $\ker b|_W =
W^{\perp}_b$.

Let now $\square$ be an affine distribution on a manifold $Q$. A form
$\omega\in \Lambda^2(Q)$ is called \emph{nondegenerate} on $\square$
if $b=\omega_x$ is nondegenerate on $U=\square_x$, for all $x\in Q$.
In such a case, there exists a unique vector field $Y \in D(Q)$ such
that $Y_x$ is the canonical displacement of $\square_x$ with respect
to~$\omega_x$, for all $x\in Q$. The vector field
$Y=Y_{\square,\omega}$ is called the \emph{canonical displacement} of
$\square$ with respect to~$\omega$.  If $\omega\in \Lambda^2(Q)$ is
nondegenerate on $\square$, then one has the isomorphism of vector
bundles
\[
\gamma = \gamma_{\Delta^0,\omega}:{\Delta^0}^*\rightarrow \Delta^0,
\qquad \gamma_x =
-(\top_{\Delta^0_x,\omega_x})^{-1}:{\Delta_x^0}^*\rightarrow
\Delta^0_x.
\]
Passing to sections of these bundles, one gets the isomorphism of
$C^{\infty}(Q)$-modules $\Gamma_{\Delta^0,\omega} : \Lambda^1_{\Delta^0}(Q)
\rightarrow D_{\Delta^0}(Q)$ defined by
\begin{equation}\label{eq:Ham-map}
  \Gamma = \Gamma_{\Delta^0,\omega}(\varrho)(x) = \gamma(\varrho(x)), \quad
  \varrho\in\Lambda^1_{\Delta^0}(Q).
\end{equation}

\begin{definition}\label{dfn:partial-symplectic-structure}
  A \emph{partial symplectic structure} on a manifold $Q$ is a pair
  $(\square , \omega)$ consisting of an affine distribution $\square$
  on $Q$ and a closed 2-form $\omega\in\Lambda^2(Q)$ which is
  nondegenerate on $\square$.
\end{definition}

Given a partial symplectic structure $\Theta=(\square,\omega)$, we
will use the subscript $\Theta$ to denote the associated objects:
$\square_{\Theta}=\square, \omega_{\Theta}=\omega,
Y_{\Theta}=Y_{\square,\omega}$ and $\Delta^0_{\Theta}$ for the
distribution canonically associated to $\square$. We also write
\[
r_{\Theta}=r_{\Delta^0_{\Theta}}, \quad
D_{\Theta}=D_{\Delta^0_{\Theta}}(Q),\quad
\Lambda^1_{\Theta}=\Lambda^1_{\Delta^0_{\Theta}}(Q), \quad
\Gamma_{\Theta}=\Gamma_{\Delta^0_{\Theta},\omega}:
\Lambda^1_{\Theta}\rightarrow D_{\Theta}.
\]
In the partial symplectic formalism, the elements of
$C^{\infty}(Q)$-modules $D_{\Theta}$ and $\Lambda^1_{\Theta}$ may be
viewed as ``constrained'' vector fields and differential forms,
respectively. The \emph{constrained Hamiltonian vector field}
associated with a Hamiltonian function $H\in C^\infty(Q)$ is defined
as
\begin{align}\label{eq:constrained-H-vector-field}
  X_H^{\Theta} = \Gamma_{\Theta}(r_{\Theta}(dH))+Y_{\Theta}.
\end{align}
The \emph{almost-Poisson} bracket associated to the partial symplectic
structure $\Theta$ is
\[
\{f,g\}_{\Theta} = \Gamma_{\Theta}(r_{\Theta}(df))
(g)=X_{f}^{\Theta}-Y_{\Theta}(g) \, , \quad f, g \in C^{\infty}(Q) \,
.
\]
The wording ``almost'' here refers to the fact that this bracket does
not satisfy in general the Jacobi identity. However, it is still
skew-symmetric and a bi-derivation.

\begin{definition}\label{dfn:partial-transversal}
  Let $\Theta=(\square, \omega)$ be a partial symplectic structure on
  a manifold $Q$. A hypersurface $B\subset Q$ is \emph{transversal} to
  $\Theta$ (or to $\square$) if the affine subspaces $T_yB$ and
  $\square_y$ of $T_yQ$ are transversal for any $y\in C$.
\end{definition}

If $B$ is \emph{transversal} to $\Theta$, then $T_yB \cap \square_y$
is of codimension $1$ in $\square _y$. If $\Theta$ is a partial
symplectic structure on $C \subset T^*M$, we shall extend this
terminology by saying that $\Theta$ is \emph{transversal} to a
hypersurface $\tilde{B}$ in $T^*M$ if $\tilde{B}$ is transversal to
$C$, so that $B=\tilde{B}\cap C$ is a hypersurface in $C$, and $B$ is
transversal to $\Theta$.


\section{Dynamics of nonholonomic Hamiltonian
  systems}\label{se:dynamics-nonholonomic}

In this section, we formulate the dynamics of nonholonomic Hamiltonian
systems using the partial symplectic formalism. We show how, under
some technical conditions, any Hamiltonian system subject to
nonholonomic constraints possesses an associated partial symplectic
structure. Then, we analyze the cases of systems with instantaneous
nonholonomic constraints, and systems exhibiting discontinuities.

\subsection{The partial symplectic structure associated with a
  constrained Hamiltonian system}\label{se:partial-ass}

Let $(M,H,C)$ be a constrained Hamiltonian system. Our goal is to
associate with it a partial symplectic structure $\Theta$ on the
``constrained'' cotangent bundle $C$ in such a way the corresponding
constrained Hamiltonian field $X_H^{\Theta}$ gives the desired
nonholonomic dynamics. With this purpose, consider the
\emph{constrained symplectic form} defined by the restriction of the
``free'' symplectic form $\Omega_M$ to $C$
\begin{align}\label{eq:partial-form}
  \omega_{\Theta} = j^*(\Omega_M) \, ,
\end{align}
with $j:C \hookrightarrow T^*M$ the canonical inclusion. The next step is
to construct a suitable affine distribution $\square_{\Theta}$ on $C$. A
natural non-singularity requirement on $C$ is asking for the regularity of
the map $\alpha_H|_C$.  This is the reason why we assume that $\alpha_H|_C$
is an immersion, i.e., that the differential $d_y\alpha_H$ is nonsingular
for any $y\in C$.  Since $\alpha_H:T^*M\rightarrow TM$ is fibered, this
assumption implies that the map $(\alpha_H)_x:C_x\rightarrow TM$ is an
immersion for any $x\in M$ and vice versa.

Let $y\in C$ and $x=\pi(y)$.  Let $\Pi_y$ be the affine subspace of
$T_xM$ tangent to $\alpha_H(C_x)$ at $z=\alpha_H(y)$. Since by the
above assumption $\alpha_H|_C$ is an immersion, $\dim\Pi_y=\dim
C_x=n-m$.  Consider the affine distribution $\square_\Theta$ on $C$
defined by
\begin{equation}\label{eq:partial-distribution}
  {\square_\Theta}_y = \{\xi\in T_yC \mid d_y \pi(\xi)\in \Pi_y \} \subset
  T_y C.
\end{equation}
Since $d_y(\pi|_C)$ is surjective, the codimension of ${\square_\Theta}_y$
in $T_yC$ is equal to the codimension of $\Pi_y$ in $T_xM$, i.e., to $m$.
Therefore, $\dim{\square_\Theta}_y=2(n-m)$.  It is not difficult to see now
that if the form $\omega_{\Theta}$ is nondegenerate on the
distribution~$\square_\Theta$, then $\alpha_H|_C$ is an immersion.

\begin{proposition}\label{prop:relation}
  Let $S = (M,H,C)$ be a constrained Hamiltonian system.  Then,
  $\alpha_H|_C$ is an immersion if the pair $(\Delta_{\Theta},
  \omega_{\Theta})$ defined by equations~\eqref{eq:partial-form}
  and~\eqref{eq:partial-distribution} is a partial symplectic structure.
\end{proposition}

The converse, however, is in general not true. Since the partial symplectic
 structure associated with $S=(M,H,C)$ is determined by $H$ and $C$, we
 will simply denote it by $\Theta(H,C)=(\square_{H,C}, \omega_C)$.

 For most Hamiltonian systems (including those coming from Mechanics),
 the anti-Legendre map $\alpha_H$ is regular not only when restricted
 to $C$, but on the whole space $T^*M$. In this is the case, and the
 Hamiltonian system comes from a Lagrangian system, one can indeed
 show that the condition of $\omega_{\Theta}$ being nondegenerate on
 the distribution~$\square_\Theta$ is equivalent to the so-called
 \emph{compatibility condition}~\cite{BaSn,LeMa}.  Therefore,
 Proposition~\ref{prop:relation} establishes a link between the
 classical partial symplectic formalism introduced in~\cite{BoVi} and
 more recent approaches as explained, for instance, in~\cite{Co}.
 Also note that the class of mechanical systems automatically verifies
 the compatibility condition, therefore admitting both formulations.
 Indeed, for mechanical systems, the conditions in
 Proposition~\ref{prop:relation} are equivalent.

\begin{definition}\label{dfn:nonholonomic-Hamiltonian-system}
  A \emph{nonholonomic Hamiltonian system} on a manifold $M$ is a
  constrained system $(M,H,C)$, $H\in C^{\infty}(T^*M)$, $C\subset T^*M$
  such that $\Theta(H,C)=(\square_{H,C}, \omega_C)$ is a partial symplectic
  structure.
\end{definition}

The dynamics of a nonholonomic Hamiltonian system is given by the
constrained Hamiltonian vector field $X^{\Theta}_H$ with respect to the
partial symplectic structure $\Theta=\Theta(H,C)$ (cf.
equation~\eqref{eq:constrained-H-vector-field}). This vector field will be
denoted by $X_{H,C}$. Under regularity of the map $\alpha_H$, $X_{H,C}$
reads in canonical coordinates
\begin{align*}
  X_{H,C} & = \pder{H}{p_a} \pder{}{q^a} - \left( \pder{H}{q^a} + {\cal
    C}_{ij} \left( \pder{H}{p_b} \pder{\Phi^j}{q^b} - \pder{H}{q^b}
    \pder{\Phi^j}{p_b} \right) \pder{\Phi^i}{p_c} {\cal H}_{ca} \right)
    \pder{}{p_a} \, ,
\end{align*}
where the matrices $({\cal H}_{ab})$ and $({\cal C}_{ij})$ are defined
by
\[
({\cal H}_{ab}) = \left( \frac{\partial^2 H}{\partial p_a \partial
    p_b} \right)^{-1} \hspace*{-8pt} , \quad ({\cal C}_{ij}) = \left(
  \pder{\Phi^i}{p_a} {\cal H}_{ab} \pder{\Phi^j}{p_b} \right)^{-1}
\hspace*{-8pt} .
\]
Observe that the \emph{force of reaction of nonholonomic constraints}
(see Section~\ref{se:dynamics}) in the partial symplectic framework is
defined \emph{a posteriori} as the difference between the constrained
and the free Hamiltonian vector fields, $X_{H,C}-X_H$.  Also, note
that the almost-Poisson bracket associated with the partial symplectic
structure $\Theta_{H,C}$ coincides with the so-called nonholonomic
bracket~\cite{Co,ScMa}.

\paragraph{Transversality}

It is convenient to adapt the terminology related to the notion of
transversality discussed in Section~\ref{se:partial} to the context of
nonholonomic Hamiltonian systems. First, we shall say that a nonholonomic
Hamiltonian system $S=(M,H,C)$ is \emph{transversal} to a hypersurface $B$
in $T^*M$ if the underlying partial symplectic structure $\Theta(H,C)$ is
transversal to $B$.  Second, if $N$ is a hypersurface in $M$, we shall say
that $S$ is \emph{transversal to $N$} if $S$ is transversal to the
hypersurface $T^*_N M$.  The following result follows from the definition
of the partial symplectic structure $\Theta(H,C)$.

\begin{proposition}\label{prop:trans}
  A nonholonomic Hamiltonian system $(M,H,C)$ is transversal to a
  hypersurface $N\subset M$ if and only if $\Pi_y \subset
  T_{\pi(y)}M$, the affine subspace of $T_{\pi(y)}M$ tangent to
  $\alpha_H(C_{\pi(y)})$ at $z=\alpha_H(y)$, is transversal to
  $T_{\pi(y)}N\subset T_{\pi(y)}M$ for all $y\in C$.
\end{proposition}

\subsection{Instantaneous partial symplectic
  structures}\label{se:inst-partial}

It is intuitive to think that when a trajectory of a Hamiltonian
system $(M,H,C)$ crosses a \emph{critical} hypersurface $N$ in the
configuration manifold $M$, its phase space reduces to $T^*N$.
Moreover, it could possibly be subject to additional instantaneous
constraints along $N$.  In the language of our approach, this idea is
naturally expressed by saying that all such \emph{critical} states
constitute a nonholonomic Hamiltonian system on $N$. Since $T^*N$ is
not naturally embedded into $T^*M$, a realization of this idea is not
completely straightforward.  What one really needs is a partial
symplectic structure on the manifold of instantaneous constraints
$C^{\text{inst}}$ which, by definition, is a submanifold of $T^*M$

Namely, let $C^{\text{inst}}$ be a set of instantaneous constraints
along $N$ imposed on $(M,H,C)$ (cf.
Definition~\ref{dfn:inst-constraints}).  Take $y\in C^{\text{inst}}$.
Let $x=\pi(y)$ and denote by $\Pi_y^{\text{inst}}$ the affine subspace
of $T_xN\subset T_xM$ tangent to $\alpha_H(C^{\text{inst}}_x)$ at
$\alpha_H(y)$. Consider the 2-form $\omega_{\Theta^{\text{inst}}}$ and
the affine distribution $\square_{\Theta^{\text{inst}}}$ on
$C^{\text{inst}}$ defined by
\begin{equation}\label{eq:inst-square}
  \omega_{\Theta^{\text{inst}}} = j^*(\Omega_M) \, , \quad
  {\square_{\Theta^{\text{inst}}}}_y = \{\xi\in T_yC^{\text{inst}}\mid
  d_y\pi(\xi)\in \Pi_y^{\text{inst}}\} \subset T_yC^{\text{inst}},
\end{equation}
with $j:C^{\text{inst}} \hookrightarrow T^*M$ the canonical inclusion.  We
 then have the following definition.

\begin{definition}
  Let $(M,H,C)$ be a nonholonomic Hamiltonian system and
  let~$C^{\text{inst}}$ be a set of instantaneous constraints along a
  hypersurface $N \subset M$.  The pair
  $(\square_{\Theta^{\text{inst}}},\omega_{\Theta^{\text{inst}}})$ defined
  by~\eqref{eq:inst-square} is called the \emph{instantaneous partial
  symplectic structure along $N$} if $\omega_{\Theta^{\text{inst}}}$ is not
  degenerate on $\square_{\Theta^{\text{inst}}}$. If this is the
  case,~$C^{\text{inst}}$ is called a \emph{regular} set of instantaneous
  constraints.
\end{definition}

Note that this structure is defined by $H$, $C^{\text{inst}}$ and $N$.
To highlight this fact, we denote $\Theta^{\text{inst}} =
\Theta^{\text{inst}}(H,C^{\text{inst}},N)$.  Accordingly, we denote by
$X_{(H,C^{\text{inst}},N)}$ the constrained Hamiltonian vector field
$X^{\Theta^{\text{inst}}}_{H^{\text{inst}}}$, with $H^{\text{inst}} =
H|_{C^{\text{inst}}}$.

In general, since $C^{\text{inst}} \subset C$ by definition, one has
that $\alpha_H|_{C^{\text{inst}}}$ is an immersion. For mechanical
systems, this implies that the 2-form $\omega_{\Theta^{\text{inst}}}$
is nondegenerate on $\square_{\Theta^{\text{inst}}}$, and therefore
any set of instantaneous constraints is regular.

In what follows, we shall only deal with regular instantaneous
nonholonomic constraints.  A natural class of instantaneous structures
arises in the following situation of particular interest. Assume that
the nonholonomic Hamiltonian system $S=(M,H,C)$ is transversal to $N$
and that $\alpha_H$ is regular. Then $\alpha_H^{-1}(TN)$ is
transversal to $C$ and, hence,
\[
C_{(N,H)} = \alpha_H^{-1}(TN)\cap C
\]
is a submanifold of codimension $2$ in $C$. Note that $C_{(N,H)}$ is a
set of instantaneous nonholonomic constraints on $S$ along $N$. By
construction, the codimension of ${\square_{N}}_y = \square_y \cap
T_y(C_{(N,H)})$ in $\square_y$ is also $2$ and $\Omega_y$ is
nondegenerate when restricted to ${\square_{N}}_y$.  Therefore, the
affine distribution $\square_N$ and the 2-form $\Omega|_N$ endow
$C_{(N,H)}$ with a partial symplectic structure, which is an
instantaneous partial symplectic structure along $N$. We call it the
\emph{trace of $S$ on $N$} and denote it by $S_{(N,H)}$. In the special
case $N=\partial M$, we call it the \emph{boundary of $S$}, and denote
it by $\partial S$, i.e., $\partial S=S_{(\partial M,H)}$. We will
denote the constrained Hamiltonian vector field with respect to the
trace (resp., boundary) as $X^{\text{tr}}=X^{\text{tr}}_{(H,C,N)}$
(resp., $X^{\partial} = X^{\partial}_{(H,C,\partial M)}$).

\subsection{Discontinuous nonholonomic
  systems}\label{se:discontinuity}

An impulsive behavior of a Hamiltonian system occurs when its
trajectory ``tries'' to go across a \emph{critical} hypersurface~$N$
in the configuration space~$M$. In such an instant, the system may be
forced to drastically change its constraints, to pass under the
control of another Hamiltonian and/or to be eventually subject to
additional instantaneous constraints. Such situations may be
interpreted as \emph{discontinuities} on both the constraints and the
Hamiltonian of the system.  Below, we formalize these concepts
properly via the notion of cutting-up.

\begin{definition}\label{dfn:cutting-up}
  Let $N\subset M$ be a hypersurface of $M$ with $N\cap \partial M =
  \emptyset$. A pair $(\hat{M},\varsigma)$, $\varsigma: \hat{M}
  \rightarrow M$, is called a \emph{cutting-up of $M$ along $N$} if
  \begin{description}
  \item [(i)] $\hat{N}=\varsigma^{-1}(N)\subset \partial\hat{M}$;
  \item [(ii)] $\varsigma$ maps $\hat{M}\setminus\hat{N}$
    diffeomorphically onto $M\setminus N$;
  \item [(iii)] $\varsigma|_{\hat{N}}:\hat{N}\rightarrow N$ is a
    double covering of $N$.
  \end{description}
\end{definition}

Note that, by definition, $\varsigma$ is a local diffeomorphism.
Cuttings-up for a given $N$ exist and are equivalent one to each
other. If $N$ divides $M$ into two parts, say,~$M_+$ and~$M_-$,
i.e.,~$M=M_+\cup M_-$, $M_+\cap M_-=N$, then $\hat{M}$ may be viewed
as the disjoint union of~$M_+$ and~$M_-$, and $\varsigma$ as the map
that matches them together along the common border $N$. Locally any
cutting-up is of this form.

For our purposes, it is important to realize that, if
$(\hat{M},\varsigma)$ is a cutting-up of $M$ along $N$, then
$(T^*\hat{M}, T^*\varsigma)$ is a cutting-up of~$T^*M$ along the
hypersurface $T^*_N M$.  Here, $T^*\varsigma$ denotes the dual of the
inverse of the isomorphism $d_z\varsigma : T_z\hat{M}\rightarrow
T_{\varsigma(z)}M$, for all $z\in \hat{M}$. In the following
definition, we introduce the class of Hamiltonian systems we shall be
dealing with throughout this paper.

\begin{definition}\label{dfn:discontinuous}
  Let $N\subset M$ be a hypersurface of $M$ with $N\cap \partial M =
  \emptyset$ and let $(\hat{M},\varsigma)$ be a cutting-up of $M$
  along $N$.  A \emph{nonholonomic Hamiltonian system discontinuous
    along $N$}, denoted $S=(M,H,C\mid N)$, is the direct image with
  respect to $T^*\varsigma$ of a nonholonomic Hamiltonian
  system~$(\hat{M},\hat{H},\hat{C})$. Such system is called
  \emph{regular} if $(\hat{M},\hat{H},\hat{C})$ is transversal to
  $\hat{N}$.

  A system of \emph{instantaneous nonholonomic constraints on $S$
    along $N$} is the direct image with respect to $T^*\varsigma$ of a
  set of instantaneous constraints~$\hat{C}^{\text{inst}}$ along
  $\hat{N}$ on the associated system
  $\hat{S}=(\hat{M},\hat{H},\hat{C}\mid \hat{N})$.  The \emph{trace of
    $S$ on $N$} is the direct image with respect to $T^*\varsigma$ of
  the trace of $\hat{S}$ along~$\hat{N}$.
\end{definition}

According to Definition~\ref{dfn:discontinuous}, $\hat{H}$ is a smooth
function on $T^*\hat{M}$ and~$\hat{C}$ is a submanifold of
$T^*\hat{M}$. Therefore, the \emph{direct image} of $\hat{H}$ along
the matching map $T^*\varsigma:T^*\hat{M}\rightarrow T^*M$ may be
viewed as a function on $T^*M$, which is 1-valued and smooth outside
of~$T^*_NM$ and 2-valued and smooth on~$T^*_NM$. We will continue to
use the notation $H$ for this function and will refer to it as a
\emph{discontinuous Hamiltonian along $N$}. Similarly, the direct
image~$C = T^*\varsigma(\hat{C})$ of~$\hat{C}$ will be referred to as
\emph{discontinuous nonholonomic constraints along $N$}.  Outside
of~$T^*_NM$, $C$ is a ``good'' smooth submanifold of $T^*M$, whose
boundary is an immersed submanifold of~$T^*_NM$.

The previous discussion becomes particularly simple when $N$ divides
$M$ into two parts, $M_+$ and $M_-$, as mentioned above. In such a
case, $T^*_NM$ also divides $T^*M$ into two parts, $T^*M_+$ and
$T^*M_-$, whose common boundary is $T^*_NM$. Then, a discontinuous
Hamiltonian $H$ along $N$ may be naturally seen as a pair of
Hamiltonians, say, $H_+$ and $H_-$, defined on $T^*M_+$ and $T^*M_-$,
respectively.  Similarly, a set of discontinuous nonholonomic
constraints along $N$ is regarded as a pair of sets of nonholonomic
constraints $C_{\pm}\subset T^*M_{\pm}$. Since $N$ always divides $M$
locally, this description constitutes a local picture of a
discontinuous nonholonomic Hamiltonian system along $N$.

We will continue to use the notation~$C^{\text{inst}}$ (resp.,
$S^{\text{tr}}$) for instantaneous nonholonomic constraints (resp.,
the trace of $S$) in the case of discontinuous nonholonomic systems.
As before, one may interpret~$C^{\text{inst}}$ as a 2-valued system of
instantaneous nonholonomic constraints along $N$. In the case when $N$
divides $M$ into two parts, we will distinguish between the two
branches using the notation~$C^{\text{inst}}_{\pm}$, and write also
$X_{(H,C^{\text{inst}}_{\pm},N)}$ (resp.,
$X^{\text{tr}}_{(H,C_{\pm},N)}$).

\begin{remark}\label{rem:boundary systems}
  {\rm The impulsive behavior of a Hamiltonian system is not
    necessarily related to some discontinuity. This type of phenomena
    occurs, for instance, each time that one of its trajectories
    ``strikes'' against the boundary $\partial M$ of the configuration
    space $M$.  Various kinds of collisions, impacts, etc, in
    mechanical systems are described in this way.  Otherwise said,
    impulsive behavior is characteristic of Hamiltonian systems
    \emph{with boundary}. Moreover, systems with boundary may be
    viewed as a ``limit'' case of discontinuous systems by dropping
    the requirement $N \cap \partial M = \emptyset$ and choosing
    $N=\partial M$, $M_-=\emptyset$, $M_+=M$. This allows a unified
    approach to both situations.}
\end{remark}

\section{The Transition Principle}\label{se:Transition-Principle}

In this section we discuss the formulation of the Transition Principle
for systems subject to nonholonomic constraints. We first introduce
the notions of focusing points, constrained characteristics, and in,
out and decisive points. The Transition Principle builds on these
elements to prescribe the behavior of the Hamiltonian system when one
or more of its ingredients undergoes a drastic change.

\subsection{Focusing points}\label{se:focusing}

The following simple linear result will be key for the subsequent
discussion.

\begin{lemma}\label{le:simple-linear}
  Let $y\in T^*M$, and let $W$ be an affine subspace in $T_y(T^*M)$
  such that $\Omega_y$ is nondegenerate on $W$ (hence, $\dim W = 2l$
  for certain $l$) and $\dim d_y\pi(W) = l$. Denote by $W^0 \subset
  T_y(T^*M)$ and ${d_y\pi(W)}^0 = d_y\pi(W^0) \subset T_{\pi(y)}M$ the
  linear subspaces associated with the affine spaces~$W$ and
  $d_y\pi(W)$, respectively. Then the affine subspaces~$W^\bullet = y
  + \ann{{d_y\pi(W)}^0}$ and~$W_\bullet=y+W^0\cap T_y(T_{\pi(y)}^*M)$
  in $T^*_{\pi(y)}M$ passing through $y$ are transversal.
\end{lemma}

\begin{proof}
  Since, by hypothesis, $\dim d_y\pi(W) = l$, one has
  \[
  \dim W^0\cap T_y(T_{\pi(y)}^*M)=l \quad \text{and} \quad \dim
  {d_y\pi(W)}^0=l \, .
  \]
  Now, the dimension of $\ann{{d_y\pi(W)}^0} \subset T_{\pi(y)}^*M$
  is~$n-l$. Moreover, $W^0\cap T_y(T_{\pi(y)}^*M)$ is transversal to
  $\ann{{d_y\pi(W)}^0}$ if one identifies the spaces $T^*_{\pi(y)}M$
  and~$T_y(T_{\pi(y)}^*M)$. The result now follows. \qed
\end{proof}

Consider now a nonholonomic Hamiltonian system $(M,H,C)$. Let $y\in
C$. Denote by $\square=\square_{(H,C)}$ be the affine distribution of
the corresponding partial symplectic structure $\Theta(H,C)$ (cf.
Section~\ref{se:partial-ass}). By
Definition~\ref{dfn:nonholonomic-Hamiltonian-system}, the affine
subspace $W=\square_y$ satisfies the assumptions of
Lemma~\ref{le:simple-linear} on $W$ (observe that $d_y\pi(W)$ is
precisely $\Pi_y$ in equation~\eqref{eq:partial-distribution}).
Therefore, the subspace $W^\bullet=\square_y^\bullet$ is well-defined
and we put
\[
K_y = K_y(H,C) = \square_y^\bullet \subset T_{\pi(y)}^*M \, .
\]
Moreover, it is not difficult to see that the subspace
$W_\bullet=(\square_y)_\bullet$ is identical to $T_yC_{\pi(y)}$.  This
shows that~$K_y$ is transversal to~$C$ at $y$, and that~$\dim K_y=m$.
The \emph{crown of the nonholonomic Hamiltonian system $(M,H,C)$} is
the map
\[
\kappa = \kappa_{H,C} : C \longrightarrow A_{m}(T^*M) \, , \quad y
\mapsto K_y \, ,
\]
where $A_k(T^*M)$ denotes the manifold whose elements are
$k$-dimensional affine submanifolds contained in the fibers of the
cotangent bundle $T^*M$. One can see that the graph of the crown
$\kappa$,
\begin{align*}
  \graph{\kappa} = \setdef{(y,v)\in C\times T^*M}{v\in K_y}
\end{align*}
is a $2n$-dimensional smooth submanifold of $C\times T^*M$. Note that
$\graph{\kappa}$ is a fiber bundle over $C$ with projection
\[
p=p_{(H,C)}: \graph{\kappa} \longrightarrow C \, , \quad (y,v)\mapsto
y \, .
\]
The fiber over $y$ of this bundle is precisely $K_y$. Since $y\in
K_y$, the map
\[
\sigma : C \longrightarrow \graph{\kappa} \, , \quad y \mapsto (y,y)
\, ,
\]
is a section of $p_{(H,C)}$. Since the fibers of $p_{(H,C)}$ are
affine spaces, the bundle $\graph{\kappa} \rightarrow C$ has a natural
vector bundle structure whose zero section is $\sigma$.  Moreover,
this vector bundle is canonically isomorphic to the normal bundle of
$C$ in $T^*M$. This is due to the fact that, for any $y\in C$, the
fiber $K_y$ is transversal to $C$ at $y$. The same argument also
guarantees that the map
\[
\Xi=\Xi_{(H,C)} : \graph{\kappa} \longrightarrow T^*M \, , \quad (y,v)
\mapsto v \, ,
\]
induces a diffeomorphism of a neighborhood of the ``zero'' section
$\sigma(C)$ in $\graph{\kappa}$ onto its image.

\begin{definition}\label{dfn:focusing}
  Let $(M,H,C)$ be a nonholonomic Hamiltonian system.  Given a point
  $u\in T^*M$, its $(H,C)$-\emph{focusing locus} $F_{(H,C)}(u)$ is the
  set of all points $y\in C$ such that $u\in K_y$. In other words,
  \[
  F_{(H,C)}(u) = p_{(H,C)} \left(\Xi_{(H,C)}^{-1}(u)\right)\subset
  C_{\pi(u)} \, .
  \]
  A point in $F_{(H,C)}(u)$ is called \emph{focusing for $u$}.
\end{definition}

Standard arguments show that $\Xi_{(H,C)}$ is regular, i.e., of maximal
rank $2n$ almost everywhere, that is, with the exception of a closed subset
without interior points. Therefore, for a generic point $u\in T^*M$, the
subset $\Xi_{(H,C)}^{-1}(u)$ is discrete, and so is $F_{(H,C)}(u)$ as
well. Note also that if $u \in C$, then $u \in F_{(H,C)}(u)$.

\begin{remark}
  {\rm Focusing points can be understood as \emph{nonintegrable}
    analogs of the notion of \emph{reducing points} considered
    in~\cite{PuVi2} in connection with the Transition Principle for
    inelastic collisions.}
\end{remark}

\begin{remark}
  {\rm It is worth noticing that the concept of a focusing point makes
    also sense in the absence of constraints. Obviously, in this case
    $F_{(H,C)}(u) = \{u\}$.  Therefore there is no need to distinguish
    between the constrained and non-constrained cases in the statement
    of the Transition Principle.}
\end{remark}

If the constraints are linear, i.e., $C=\Upsilon+C^0$ with $C^0$ a linear
codistribution and $\Upsilon \in \Lambda^1(M)$ (the displacement form),
then for each $y\in T^*M$,
\begin{align*}
  T_{\pi(y)}^*M = C_{\pi(y)}^o \oplus \ann{d_y \pi(\Delta_y)} \, ,
\end{align*}
where $\Delta=\Delta_{(H,C)}$. Denote the corresponding projectors by $\P:
T_{\pi(y)}^*M \rightarrow C_{\pi(y)}^o$ and $\Q: T_{\pi(y)}^*M \rightarrow
\ann{d_y \pi (\Delta_y)}$.  Given $u \in T^*M$, one has that $z \in
F_{(H,C)}(u)$ if and only if $z\in C$ and~$\P(z)=\P(u)$.  Since $z = \P(z)
+ \Q(z) = \P(u) +\Q(\Upsilon_y)$, one has the following result.

\begin{proposition}\label{prop:focusing-affine}
  Let $(M,H,C)$ be a nonholonomic Hamiltonian system with linear
  constraints. Then, for $u \in T_y^*M$, there is a unique focusing point
  given by $F_{(H,C)}(u) = \{ \P(u) + \Q(\Upsilon_y) \}$.
\end{proposition}

\subsection{Instantaneous focusing points}\label{se:inst-focusing}

We will also need an instantaneous version of the notion of a focusing
point introduced in the previous section. For this purpose, it is
sufficient to apply the above construction to instantaneous
constraints instead of to the ``usual'' ones.  Namely, let
$C^{\text{inst}}$ be a system of regular instantaneous constraints
along $N$ (see Section~\ref{se:inst-constraints}) and
$\square^{\text{inst}} = \square_{\Theta^\text{inst}}$ be the
corresponding affine distribution (see Section~\ref{se:inst-partial}).
Following the same reasoning as above, the affine subspace $W =
\square^{\text{inst}}_y \subset T_y(T^*M)$ satisfies the assumptions
of Lemma~\ref{le:simple-linear}.  Therefore, the affine subspace
$K_y^{\text{inst}} = (\square^{\text{inst}})^\bullet$ of $T_y(T^*M)$
is well-defined, and we have all the ingredients to define the notion
of instantaneous crown and instantaneous focusing point of a system
subject to instantaneous nonholonomic constraints. For completeness,
we state the latter.

\begin{definition}\label{dfn:focusing.instantaneous}
  Let $(M,H,C)$ be a nonholonomic Hamiltonian system and let
  $C^{\text{inst}}$ be a set of instantaneous constraints along a
  hypersurface $N \subset M$. Given a point $u \in T_N^*M$, its
  $(H,C^{\text{inst}},N)$-\emph{instantaneous focusing locus}
  $F_{(H,C^{\text{inst}},N)}(u)$ is the set of all points $y\in
  C^{\text{inst}}$ such that $u \in K_y^{\text{inst}}$. In other
  words,
  \[
  F_{(H,C^{\text{inst}},N)}(u) = p_{(H,C^{\text{inst}})}
  \left(\Xi_{(H,C^{\text{inst}},N)}^{-1}(u)\right)\subset
  C^{\text{inst}}_{\pi(u)} \, .
  \]
  A point in $F_{(H,C^{\text{inst}},N)}(u)$ is called
  \emph{instantaneous focusing for $u$}.
\end{definition}

As before, if the instantaneous nonholonomic constraints are linear
$C^{\text{inst}} = \Upsilon^\text{inst}+{C^{\text{inst}}}^o$, then for each
$y\in T_N^*M$,
\begin{align*}
  T_{\pi(y)}^*M = C_{\pi(y)}^o \oplus \ann{d_y
    \pi(\Delta^{\text{inst}}_y)} \, ,
\end{align*}
where $\Delta^{\text{inst}} = \Delta_{(H,C^{\text{inst}},N)}$.
Denoting the corresponding projectors by $\P^{\text{inst}}:
T_{\pi(y)}^*M \rightarrow {C^{\text{inst}}}^o_{\pi(y)}$ and
$\Q^{\text{inst}}: T_{\pi(y)}^*M \rightarrow \ann{d_y \pi
  (\Delta^{\text{inst}}_y)}$, one has the following result.
\begin{proposition}\label{prop:instantaneous-focusing}
  Let $(M,H,C)$ be a nonholonomic Hamiltonian system and let
  $C^{\text{inst}}$ be a set of instantaneous affine constraints along a
  hypersurface $N \subset M$.  Then, for $u \in T_y^*M$, there is a unique
  instantaneous focusing point given by $ F_{(H,C^{\text{inst}},N)}(u) = \{
  \P^{\text{inst}}(u) + \Q^{\text{inst}} (\Upsilon^\text{inst}_y) \}$.
\end{proposition}


\subsection{Constrained characteristics}\label{se:characteristics}

Consider then a partial symplectic structure $\Theta=(\square,\omega)$
on a manifold $C$ which is transversal to a hypersurface $B\subset C$
(cf. Definition~\ref{dfn:partial-transversal}).  Let $\Delta^0$ denote
the linear distribution associated with $\square$.  For each $y \in
B$, consider the linear space $V=\Delta^0_y$, the hyperplane $W =
\Delta^0_y\cap T_yB$ of $V$ and the nondegenerate skew-symmetric form
$b = \omega_y|_{\Delta^0_{y}}$.  The \emph{characteristic direction at
  $y \in B$} is defined as
\[
l_y = l_y(\Theta, B) = \ker b|_W = \ker(\omega_y|_{\Delta^0_y \cap
  T_yB}) \subset \Delta^0_y \cap T_yB \, .
\]
The proof of the following result is straightforward.
\begin{lemma}
  Given a partial symplectic structure $\Theta=(\square,\omega)$ on a
  manifold $C$ and a hypersurface $B \subset C$ transversal to it, the
  distribution $y\mapsto l_y(\Theta,B)$ is one-dimensional.
\end{lemma}

\begin{definition}
  Given a partial symplectic structure $\Theta=(\square,\omega)$ on $C$ and
  a hypersurface $B \subset C$ transversal to it, $y\mapsto l_y(\Theta,B)$
  is called the \emph{characteristic distribution with respect to
  $(\Theta,B)$}, and its integral curves, denoted by $\zeta$, are the
  \emph{$(\Theta,B)$-characteristics}.
\end{definition}

We are particularly interested in the case when we have a nonholonomic
Hamiltonian system $S=(M,H,C)$, the partial symplectic structure $\Theta$
is $\Theta_{H,C}$, $N$ is a hypersurface in $M$ and $\tilde{B}=T^*_NM$,
$B=T^*_NM \cap C$. We will use the terminology $(S|N)$- or
$(H,C|N)$-\emph{characteristic} as a substitute for
$(\Theta,B)$-characteristic.  It should be emphasized that
$(H,C|N)$-characteristics are only defined when $S$ is transversal to $N$
(see Section~\ref{se:partial-ass}).

In the absence of constraints, i.e., when $(C=T^*M,\omega=\Omega)$ is
a symplectic manifold, and $\Delta$ is the trivial distribution $y
\mapsto T_yC$ on $C$, the characteristic curves are precisely the
characteristics introduced in~\cite{BoVi}. We will refer to
\emph{non-constrained characteristics} and \emph{constrained
  characteristics} when it is necessary to distinguish between the
unconstrained and the constrained cases.

\begin{remark}
  {\rm Just as non-constrained characteristics play a key role in
    describing holonomic elastic collisions, and reflection and
    refraction phenomena of rays of light~\cite{BoVi,PuVi}, the
    constrained characteristics will be fundamental in describing the
    ``elastic part'' of nonholonomic impulsive phenomena. What is
    meant by ``elastic part'' will become clear in
    Section~\ref{se:decisive} when describing decisive points.}
\end{remark}

If the constraints $C$ are affine, then the $(H,C|N)$-characteristic
passing through a point $y \in C$, $\pi(y)\in N$, is described in a
particularly simple way. Namely, following
Proposition~\ref{prop:trans}, it is not difficult to see that the
$(H,C|N)$-characteristic passing through $y$ is given by
\[
\zeta_y = y+C^0_{\pi (y)} \cap \ann{d_y \alpha_H (C^0) \cap
  T_{\pi(y)}N} \, ,
\]
with $C^0$ being the linear codistribution associated to $C$. In
particular, in the absence of constraints, $C=T^*M$ and the
characteristics are straight lines in $T_x^*M$ parallel to
$\ann{T_xN}$, $x\in N$.


\subsection{In, out and trapping points}\label{se:in-out-trapping}

Here, we first introduce some concepts concerning the behavior of a
vector field in a neighborhood of the boundary of its supporting
manifold. We then discuss the notions in, out and trapping points.

Let $Q$ be a manifold with boundary and $X$ a vector field on $Q$. A
point $y\in \partial Q$ is called a \emph{$j$th order in point for
  $X$} if there exists a trajectory of $X$, $\beta :[0,a]\rightarrow
Q$, $a>0$ such that
\begin{align*}
 y=\beta (0) \quad \text{and} \quad \beta (t) \notin \partial Q \, , \quad
 \text{for} \; 0<t\leq a \, ,
\end{align*}
and $\beta$ is $j$th order tangent to $\partial Q$ at $y$.  A
\emph{$j$th order out point for $X$} is a $j$th order in point for
$-X$. In the dynamical context we have in mind, in and out points of
$0$th order are the most important. It is easy to see that $y\in
\partial Q$ is a $0$th order in point (resp., out point) for $X$ if
the vector $X_y$ is transversal to $\partial Q$ and directed inside
(resp., outside) of $Q$. A point that lies on a trajectory of $X$
which is entirely contained in $\partial Q$ is called a \emph{trapping
  point for $X$}.

Let $\partial Q^j = \partial Q^j (X)$ denote the subset of all points
of $\partial Q$ where $X$ is $j$th order tangent to $\partial Q$, and
$\partial Q_>^j=\partial Q_>^j(X)$ (resp., $\partial Q^j_<=\partial
Q_<^j (X)$) the set of all $j$th order in points (resp., out points)
for $X$.  Note that $\partial Q^j \supset \partial Q^{j+1}$ and
\begin{equation}\label{eq:DDD}
  \partial Q^j\setminus(\partial Q_>^j\cup \partial Q_<^j) \subset
  \partial Q^{j+1}.
\end{equation}
In a generic situation, $\partial Q^j$ is a submanifold (with
singularities) of codimension $j$ in $\partial Q$, which is divided by
$\partial Q^{j+1}$ into two parts, $\partial Q_>^j \setminus \partial
Q_>^{j+1}$ and $\partial Q_<^j \setminus \partial Q_<^{j+1}$.  An
analytical description of the previous discussion is obtained by
choosing a smooth function $f$ on $Q$ with $f\geq 0$ and $d_zf\neq 0$,
for all $z\in \partial Q$ such that $\partial Q=\{f=0\}$ (which always
exists locally).  Then
\begin{align*}
  \partial Q^j & = \setdef{z\in Q}{f(z)=0,\; X(f)(z)=0,\dots,
    X^j(f)(z) = 0} \, , \\
  \partial Q_>^j\setminus \partial Q_>^{j+1} &= \partial Q^j \cap
  \setdef{z \in Q}{X^{j+1}(f)(z)>0} \, ,\\
  \partial Q_<^j\setminus \partial Q_<^{j+1} &= \partial Q^j \cap
  \setdef{z \in Q}{X^{j+1}(f)(z)<0} \, .
\end{align*}
The vector field $X$ is said to be \emph{regular with respect to
  $\partial Q$} when the inclusion in equation~\eqref{eq:DDD} is an
equality for all $j\geq 0$. This is a generic property of vector
fields. In such a case, the chain of inclusions
\[
\partial Q = \partial Q^0 \supset \partial Q^1 \supset \dots \supset
\partial Q^j \supset \dots \supset \partial Q^n
\]
is a \emph{stratification} of $\partial Q$ whose strata are $\partial
Q_>^j\setminus \partial Q_>^{j+1}$ and $\partial Q_<^j\setminus
\partial Q_<^{j+1}$. Note also that the set of trapping points
precisely corresponds to $\partial Q^n$.

Consider now a discontinuous nonholonomic Hamiltonian system
$(M,H,C\mid N)$ and the corresponding cutting-up
$(\hat{M},\varsigma)$. Then we can resort to the previous discussion
with the manifold $Q = \hat{C} \subset T^*\hat{M}$ and the vector
field $X=X_{\hat{H},\hat{C}}$. Recall that $\hat{N} \subset \partial
\hat{M}$ and $\partial\hat{C}=\hat{C}\cap T^*_{\partial \hat{M}}
\hat{M}$.

\begin{definition}\label{dfn:in-out}
  Let $S=(M,H,C\mid N)$ be a discontinuous nonholonomic system and
  denote by $(\hat{M},\varsigma)$ the associated cutting-up.  A point
  $y\in T^*_NM$ is called an \emph{in point} (resp., an \emph{out
    point}) of $S$ if there exists $z \in T_{\hat{N}}\hat{M}$ such
  that $y=\varsigma (z)$ and $z$ is an in point (resp., an out point)
  of $X_{\hat{H},\hat{C}}$ with respect to $\partial\hat{C}$.
\end{definition}

By definition, the map $T^*\varsigma$ restricted to $\partial\hat{C}$
is an immersion. A point in $T_N^*M$ may turn out to be an in and an
out point at the same time.  To resolve this ambiguity, the branch of
$T^*\varsigma$ to which such a point belongs must be taken into
consideration. This distinction is easily described in the case when
$N$ divides $M$ into two parts. In fact, in this case the system
$(M,H,C\mid N)$ may be viewed as a couple of nonholonomic Hamiltonian
systems $(M_{\pm},H_{\pm}, C_{\pm})$, with the common boundary
$\partial M_{\pm}=N$, and where $H_{\pm}\in C^{\infty}(M_{\pm})$ and
$C_{\pm}\subset M_{\pm}$ (cf. Section~\ref{se:discontinuity}). An in
(resp., out, or trapping) point of the vector field $X_{H_+,C_{+}}$
with respect to the boundary $\partial C_{+}$ is called an
\emph{plus-in} (resp., \emph{plus-out}, or \emph{plus-trapping)}
\emph{point}. Analogous definitions are established for $\varepsilon =
-$. In this way, the notions of plus-in point, minus-in point, etc,
introduced in~\cite{BoVi} for the unconstrained situation are
generalized to the constrained case.  Finally, we observe that $N$
always divides $M$ locally, and therefore the previous discussion is
always valid locally.

\subsection{Decisive points}\label{se:decisive}

At this point, we are ready to introduce the key notion of decisive
point corresponding to an out point. The construction of decisive
points depends on two elements: first, the mode (elastic or inelastic)
in which the system passes through the critical state and, second, the
continuity and differentiability properties of the Hamiltonian. Below,
we will limit our discussion to the two most relevant situations, just
to avoid not very instructive technicalities arising in the most
general context.  The first one is the case when the Hamiltonian is
smooth and only the constraints are discontinuous along the critical
hypersurface. The second one concerns discontinuous Hamiltonians and
not necessarily discontinuous constraints. It is worth stressing that
the first situation can not be considered as a particular case of the
second one, i.e., that the notion of a decisive point is not
``continuous'' in this sense. In what follows, $\varepsilon \in
\{+,-\}$ and $\bar{\varepsilon}$ stands for the opposite sign to
$\varepsilon$. Throughout the section, instantaneous constraints are
assumed to be regular.

\subsubsection{Elastic mode: change of constraints}

Here, we deal with a discontinuous nonholonomic system $(M,H,C|N)$,
where the Hamiltonian function is smooth, $H\in C^{\infty}(M)$.

\begin{definition}[Decisive points for smooth Hamiltonians and
  discontinuous constraints] \label{dfn:decisive-change-constraints}
  Let $(M,H,C|N)$ be a regular discontinuous nonholonomic system,
  with $H\in C^{\infty}(M)$ and consider a set of instantaneous
  constraints $C^{\text{inst}}$ along $N$.  Let $y$ be an
  $\varepsilon$-out point of the system. A sequence
  $(y_i,\varepsilon_i)$, $i=0,1,\dots,k$, with $y_i\in C\cap T^*_NM$
  is called $(y,\varepsilon)$-\emph{admissible} if it verifies the
  following conditions:
  \begin{description}
  \item[(i)] $(y_0,\varepsilon_0)=(y,\varepsilon)$;
  \item[(ii)] for all $i<k$, $y_{i+1}$ is a focusing point for $y_i$
    with respect to either $C_{\varepsilon_{i+1}}^{\text{inst}}$ or,
    if instantaneous constraints are absent, $C_{\varepsilon_{i+1}}$;
  \item[(iii)] $y_i$ is an $\varepsilon_i$-out point for all $i<k$ and
    $y_k$ is either an $\varepsilon_k$-in point or an
    $\varepsilon_k$-trapping point;
  \item[(iv)] the sequence of signs $\{\varepsilon_i\}$ alternates, i.e.,
    $\varepsilon_{i+1}=\bar{\varepsilon}_i$.
  \end{description}
  The \emph{end} point of an $(y,\varepsilon)$-admissible sequence,
  $(y_k,\varepsilon_k)$, is called $(y,\varepsilon)$-\emph{decisive}
  and the constrained Hamiltonian vector field $X_{H,
    C_{\varepsilon_k}}$ is referred to as the vector field
  \emph{corresponding} to it.
\end{definition}

\begin{remark}\label{rem:chattering}
  {\rm The above formal description of decisive points is equivalent
    to the following iterative procedure. Take, for instance, a
    plus-out point $y$.  Then, according to
    Definition~\ref{dfn:decisive-change-constraints}, all focusing
    with respect to $C_{\bar{\varepsilon}}^{\text{inst}}$ (resp., to
    $C_{\bar{\varepsilon}}$) minus-in and minus-trapping points are
    decisive. On the other hand, the procedure continues by restarting
    from any of the remaining focusing points that are minus-out
    points, and so on.  In some situations, this process may turn out
    to be infinite. At the present time, however, it is not clear
    whether that kind of phenomena can occur, say in propagation of
    singularities or similar processes.}
\end{remark}

\subsubsection{Elastic mode: discontinuous  Hamiltonians}

In this case, decisive points are constructed on the basis of an
iterative procedure whose single steps are either of \emph{reflective}
or of \emph{refractive} type, as described below. Consider a regular
discontinuous nonholonomic system $S=(M,H,C|N)$, which might be
subject to additional instantaneous constraints $C^{\text{inst}}$
along $N$. Let $y \in C \cap T_N^*M$ be an $\varepsilon$-out point.

\paragraph{Reflective step}
\begin{description}
\item [1-st move:] $y\Rightarrow z$, where $z$ is a point in the
  constrained characteristic
  $\zeta_y(H_{\varepsilon},C_{\varepsilon})$ such that
  $H_{\varepsilon}(z)=H_{\varepsilon}(y)$.
\item [2-nd move:] $z\Rightarrow u$, where $u$ is a focusing point for
  $z$ with respect to either $C^{\text{inst}}_{\varepsilon}$ or, if
  $\varepsilon$-instantaneous constraints are absent,
  $C_{\varepsilon}$.
\end{description}

\paragraph{Refractive step}
\begin{description}
\item [1-st move:] $y\Rightarrow z$, where $z$ is a point of the
  constrained characteristic $\zeta_y(H_{{\varepsilon}}
  ,C_{{\varepsilon}})$ and such that
  $H_{\bar{\varepsilon}}(z)=H_{\varepsilon}(y)$.
\item [2-nd move:] $z\Rightarrow u$, where $u$ is a focusing point for
  $z$ with respect to either $C^{\text{inst}}_{\bar{\varepsilon}}$ or,
  if $\bar{\varepsilon}$-instantaneous constraints are absent,
  $C_{\bar{\varepsilon}}$.
\end{description}

With a slight abuse of language, we shall say that $(y,\varepsilon)$
is the \emph{initial point of the step} and $(u,\varepsilon)$ (resp.,
$(u,\bar{\varepsilon})$) is the \emph{end point of the step} if the
scenario is reflective (resp., refractive).

\begin{definition}[Decisive points for discontinuous
  Hamiltonians]\label{dfn:decisive-dis-elastic} Consider a regular
  discontinuous nonholonomic system $(M,H,C|N)$. Let $C^{\text{inst}}$ be a
  set of instantaneous constraints along $N$.  Let $y$ be an
  $\varepsilon$-out point. A sequence $(y_i,\varepsilon_i)$,
  $i=0,1,\dots,k$, is called $(y,\varepsilon)$-\emph{admissible} if
  \begin{description}
  \item[(i)] $(y_0,\varepsilon_0)=(y,\varepsilon)$;
  \item[(ii)] $(y_i,\varepsilon_i)$ and $(y_{i+1},\varepsilon_{i+1})$
    are the initial and the end points of a step, respectively;
  \item[(iii)] $y_i$ is an $\varepsilon_i$-out point, $0\leq i<k$, and
    $y_k$ is an $\varepsilon_k$-in point or an
    $\varepsilon_k$-trapping point.
  \end{description}
  The end point $(y_k, \varepsilon_k)$ of an admissible sequence is
  called $(y,\varepsilon)$-decisive and the constrained Hamiltonian
  vector field $X_{H, C_{\varepsilon_k}}$ is referred to as the vector
  field \emph{corresponding} to it.
\end{definition}

If the Hamiltonian is discontinuous and the constraints are linear, i.e.,
$C\subset T^*M$ is a smooth linear submanifold, and the instantaneous
constraints are absent, the previous definition of decisive points becomes
much simpler, as the following result shows.

\begin{proposition}\label{prop:decisive-change-Hamiltonian}%
  Let $(M,H,C|N)$ be a regular discontinuous nonholonomic Hamiltonian
  system with smooth linear constraints.  Let $y$ be an $\varepsilon$-out
  point. The $(y,\varepsilon )$-\emph{decisive} points are the in and the
  trapping points belonging to the intersection of the constrained
  characteristic $\zeta_y$ passing through $y$ with the set $\setdef{z\in
  C}{H_{\pm}(z)=H_{\varepsilon}(y)}$.
\end{proposition}

\begin{proof}
  Let $y$ be an $\varepsilon$-out point and denote by $\{z_1,\dots,z_s\}$
  (resp, $\{\ov{z}_1,\dots,\ov{z}_{\bar{s}}\}$) the points belonging to the
  intersection of the constrained characteristic $\zeta_y$ passing through
  $y$ with the set $\setdef{z\in C}{H_{\varepsilon}(z)=H_{\varepsilon}(y)}$
  (resp. with $\setdef{z\in C}{H_{\bar{\varepsilon}}(z) =
  H_{\varepsilon}(y)}$).  Since the constraints are smooth, then $u=z$ in
  the $2$nd-move of both a reflective and a refractive step. Now, for any
  $j \in \until{s}$, the intersection of the constrained characteristic
  passing $\zeta_{z_j}$ through $z_j$ with the set $\setdef{z\in
  C}{H_{\varepsilon}(z)=H_{\varepsilon}(z_j)}$ (resp.  with $\setdef{z\in
  C}{H_{\bar{\varepsilon}}(z)=H_{\varepsilon}(y)}$) is again
  $\{z_1,\dots,z_s\}$ (resp, $\{\ov{z}_1,\dots,\ov{z}_{\bar{s}}\}$). The
  same observation holds for any $\ov{z}_j$, $j \in \until{\ov{s}}$.  The
  result now follows from Definition~\ref{dfn:decisive-dis-elastic}.  \qed
\end{proof}

\begin{remark}\label{rem:decisive-boundary}
  {\rm The introduced terminology remains valid for nonholonomic
    systems with boundary (cf. Remark~\ref{rem:boundary systems}). In
    such a case, one has to formally put
    \[
    M_- = \emptyset \, , \quad M_+=M \, , \quad N=\partial M \, ,
    \quad H_-=\infty\, , \quad H_+ = H \, .
    \]
    This type of geometric data occurs in describing various collision
    phenomena.  }
\end{remark}

\subsubsection{Inelastic mode: change of constraints}

As in the elastic case, we first deal with the case when the
Hamiltonian $H$ is smooth. We treat an inelastic behavior of the
system as the passage under the control of either the instantaneous
discontinuous nonholonomic system or, if instantaneous constraints are
absent, the discontinuous boundary system.  In this and subsequent
sections, the following shorthand notation will be used (cf.
Sections~\ref{se:inst-partial} and~\ref{se:discontinuity})
\[
\begin{array}{rlllrll}
  C^{\text{inst,tr}}_{\varepsilon} &=& C^{\text{inst}}_\varepsilon \cap
  \alpha_{H_\varepsilon}^{-1} (TN) \, , &\quad&
  X^{\text{inst,tr}}_{\varepsilon} &=&
  X_{(H_{\varepsilon},C^{\text{inst,tr}}_{\varepsilon},N)} \,
  , \\
  C^{\text{tr}}_{\varepsilon} &=& {C_\varepsilon}_{(N,H_{\varepsilon})}
  \, , &\quad& X^{\text{tr}}_{\varepsilon} &=&
  X^{\text{tr}}_{(H_{\varepsilon},C_{\varepsilon},N)} \, .
\end{array}
\]
We also use this notation when the Hamiltonian~$H$ is smooth, i.e.,
$H_{\pm} = H$.

\begin{definition}[Decisive points for smooth Hamiltonians and
  discontinuous constraints] \label{dfn:decisive-change-inelastic}
  Consider a regular discontinuous nonholonomic system $(M,H,C|N)$.
  Let $C^{\text{inst}}$ be a set of instantaneous constraints along
  $N$. Let $y$ be an $\varepsilon$-out point. An
  $(y,\varepsilon)$-\emph{decisive point} is a focusing point for $y$
  with respect to either $C_{\bar{\varepsilon}}^{\text{inst,tr}}$ or,
  if the instantaneous constraints are absent,
  $C^{\text{tr}}_{\bar{\varepsilon}}$.  The constrained Hamiltonian
  vector field $X^{\text{inst,tr}}_{\bar{\varepsilon}}$, respectively,
  $X^{\text{tr}}_{\bar{\varepsilon}}$ is referred to as the
  \emph{corresponding} vector field.
\end{definition}

\subsubsection{Inelastic mode: discontinuous Hamiltonians}

As in the elastic case, decisive points are constructed on the basis
of \emph{reflective} or \emph{refractive} steps, as we now describe.

\paragraph{Reflected falling step}
\begin{description}
\item [1-st move:] $y\Rightarrow z$, where $z$ is a point of the
  constrained characteristic $\zeta_y(H_{\varepsilon},
  C_{\varepsilon})$ such that $H_{\varepsilon}(z) =
  H_{\varepsilon}(y)$.
\item [2-nd move:] $z\Rightarrow u$, where $u$ is a focusing point for
  $z$ with respect to either $C^{\text{inst,tr}}_{\varepsilon}$ or, if
  $\varepsilon$-instantaneous constraints are absent, $
  C^{\text{tr}}_{\varepsilon}$.
\end{description}

\paragraph{Refracted falling step}
\begin{description}
\item [1-st move:] $y\Rightarrow z$, where $z$ is a point of the
  constrained characteristic $\zeta_y(H_{{\varepsilon}} ,
  C_{{\varepsilon}})$ such that $H_{\bar{\varepsilon}}(z) =
  H_{\varepsilon}(y)$.
\item [2-nd move:] $z\Rightarrow u$, where $u$ is a focusing point for
  $z$ with respect to $C^{\text{inst,tr}}_{\bar{\varepsilon}}$ or, if
  $\bar{\varepsilon}$-instantaneous constraints are absent,
  $C^{\text{tr}}_{\bar{\varepsilon}}$.
\end{description}

We shall refer to $(u,\varepsilon)$ (resp., $(u,\bar{\varepsilon})$)
as a \emph{reflected} (resp. \emph{refracted}) \emph{falling point}.

\begin{definition}[Decisive points for discontinuous Hamiltonians]
  \label{dfn:decisive-dis-inelastic}
  Consider a regular discontinuous nonholonomic system $(M,H,C|N)$.
  Let $C^{\text{inst}}$ be a set of instantaneous constraints along
  $N$.  Let $y$ be an $\varepsilon$-out point. An
  $(y,\varepsilon)$-\emph{decisive} point is a falling point for $y$.
  The vector field $X^{\text{inst,tr}}_{\varepsilon}$ (resp.,
  $X^{\text{tr}}_{\varepsilon}$ if $C^{\text{inst}}_{\varepsilon} =
  \emptyset$) is called the vector field \emph{corresponding} to a
  reflected falling point.  The vector field
  $X^{\text{inst,tr}}_{\bar{\varepsilon}}$ (resp.,
  $X^{\text{tr}}_{\bar{\varepsilon}}$ if
  $C^{\text{inst}}_{\bar{\varepsilon}} = \emptyset$) is called the
  vector field \emph{corresponding} to a refracted falling point.
\end{definition}

\subsection{Transition Principle}

From a physical point of view, the Transition Principle formulated
below is an explicit description of the discontinuity of a trajectory
of a regular nonholonomic Hamiltonian system $S$ that occurs when it
traverses a \emph{critical state}.  Such a discontinuity is
interpreted as an impact, collision, reflection, refraction, etc,
depending on the physical situation modeled by the system~$S$. From a
mathematical point of view, the Transition Principle corresponds to
the definition of the trajectory of a regular discontinuous
nonholonomic Hamiltonian system.

The elastic or inelastic character of the impulsive motions of an
specific physical system must be taken into account when defining the
trajectories.  Accordingly, there are two different versions of the
Transition Principle that distinguish between the two situations.  Let
$S=(M,H,C|N)$ stand for a regular discontinuous nonholonomic system
and let $C^{\text{inst}}$ be eventual instantaneous constraints
imposed on $S$ along $N$. Let $(\hat{M},\varsigma)$ be the associated
cutting-up of $M$ along $N$ (cf. Section~\ref{se:discontinuity}). The
\emph{regular part of a trajectory of the system
  $\hat{S}=(\hat{M},\hat{H},\hat{C})$} is the part of the trajectory
of the Hamiltonian vector field $X_{\hat{H},\hat{C}}$ that lies
outside $\partial\hat{M}$. The \emph{regular part of a trajectory of
  $S$} is the image by $\varsigma$ of the regular part of the
corresponding trajectory of~$\hat{S}$. At least locally, the regular
part may be viewed as a piece of the trajectory of the vector field
$X_{H_{\varepsilon},C_{\varepsilon}}$ that lies outside the
hypersurface $T^*_NM$.

\smallskip

\noindent {\bf Transition Principle.} \emph{Let
  $S=(M,H,C|N)$ be a regular discontinuous nonholonomic system and let
  $C^{\text{inst}}$ be eventual instantaneous constraints on $S$ along
  $N$. If a regular trajectory of the vector field
  $X_{H_{\varepsilon},C_{\varepsilon}}$, $\varepsilon=\pm$ reaches the
  critical hypersurface $T_N M$ at a point $y$, it then continues its
  motion from any $(y,\varepsilon)$-decisive point according to the
  chosen mode, elastic or inelastic, under the control of the
  corresponding constrained Hamiltonian vector field.}

\smallskip

Some features of the Transition Principle are worth mentioning.  First
of all, it prescribes a \emph{splitting of the trajectory} when the
number of decisive points is greater than one.  Of course, it is
difficult to imagine that a true mechanical system ``goes into
pieces'' when reaching the critical hypersurface. But it may perfectly
happen when a Hamiltonian system describes the propagation of
singularities in a fields or a continuum media. A classical example
one finds in geometrical optics when a light ray passing from one
optic medium to another splits into reflected and refracted rays (see,
for instance,~\cite{PuVi}). The trajectory may also be \emph{trapped}
by the critical hypersurface.  This happens when an ``impact'' state
$y$ possesses no $y$-decisive points.


\section{Mechanical systems}\label{se:mechanical-systems}

In this section, we particularize the previous discussion to
mechanical systems subject to affine constraints.  Let $g$ be a
Riemannian metric on $M$ and $V \in C^{\infty}(M)$, and consider the
mechanical system whose kinetic energy and potential function are
$T(q,v)=\frac{1}{2}g(v,v)$ and $V$, respectively. The corresponding
Lagrangian function is $L(q,v)=T(q,v)-V(q)$ and the Hamiltonian one is
\begin{equation}\label{eq:Hamiltonian-mechanical}
  H(q,p)=\hat{T} (q,p) + V(q) \, ,
\end{equation}
where $\hat{T}(q,p) = \frac{1}{2} \G (p,p)$, and $\G$ is the
co-metric, i.e., the metric on the cotangent bundle induced by $g$. In
a local chart $q^a$ on $M$, the local expressions of $g$ and $\G$ are
\[
g=g_{ab} dq^a \otimes dq^b \, , \quad \G = g^{ab} \pder{}{q^a} \otimes
\pder{}{q^b} \, .
\]
In the mechanical case, the Legendre transform $\Legendre_L:TM \lo
T^*M$ is a linear bundle mapping whose local description is
$\Legendre_L(q^a,\dot{q}^a) = (q^a,g_{ab}\dot{q}^b)$.

Consider an affine distribution $C = C^0 + Y$ in $T^*M$
determining some nonholonomic constraints on the system $(M,H)$.
The linearity of $\alpha_H = \Legendre_L^{-1}$ implies that the
space $\alpha_H(C) = \alpha_H(C^0) + \alpha_H(Y)$ is a
distribution of affine spaces on $M$, or otherwise said, that
$\alpha_H(C^0)$ is a linear distribution on~$M$. Throughout this
section, we will often resort to the shorthand notation $\D=
\alpha_H(C^0)$ and $\Upsilon = \alpha_H(Y)$. Now, it is easy to
verify that $T^*_q M = C^0_q \oplus \ann{\D}_q$, with associated
projectors
\[
\P_q : T^*_qM \lo C^0_q \; , \quad \Q_q : T^*_qM \lo \ann{\D}_q \, ,
\quad q \in M \, .
\]
Let $\mu_1=\mu_{1a}dq^a,\dots,\mu_m = \mu_{ma}dq^a$ be 1-forms
such that (locally) $\ann{\D} = \spn \{ \mu_1,\dots,\mu_m \}$.
Define the local function $\mu_{i0}:M \rightarrow \real$ by
$\mu_{i0} (q) = - \mu_i (\Upsilon (q))$. Then $\alpha_H(C)$ is
locally defined by the equations
\[
\mu_{ia}(q) \dot{q}^a + \mu_{i0}(q) = 0 \, , \quad 1 \le i \le m \, .
\]
Now, consider the matrices
\begin{equation}\label{eq:matrices}
  G = \left( g_{ab} \right) \, , \quad J=(\mu_{ia}) \, , \quad
  \mathcal{B} = JG^{-1}J^t \, .
\end{equation}
From the discussion after Proposition~\ref{prop:relation}, recall that
$(\Delta_{H,C},\omega_C)$ is a partial symplectic structure if and
only if $\alpha_H|_C$ is an immersion, or, equivalently, if the
compatibility condition is verified. Following~\cite{CaRa}, the latter
is equivalent to the matrix $\mathcal{B}$ being invertible.  A direct
computation give the following local expression for the projectors
$\P$ and $\Q$,
\[
\P (x) = x - \Q (x) \; , \quad \Q (x) = J^t \mathcal{B}^{-1} J G^{-1}
x \, , \quad x \in T^*M \, .
\]
Finally, let $N \subset M$ be a hypersurface and assume that the
nonholonomic system $S=(M,H,C)$ is transversal to $N$. Consider also a
set of instantaneous nonholonomic linear constraints $C^{\text{inst}}
= ({C^{\text{inst}}}^o,\Upsilon_{C^{\text{inst}}})$ imposed on $S$
along $N$. Note that $T^*_q M = {C^{\text{inst}}_q}^o \oplus
\ann{\D^{\text{inst}}}_q$, with associated projectors
\[
\P^{\text{inst}}_q : T^*_qM \lo {C^{\text{inst}}_q}^o \; , \quad
\Q^{\text{inst}}_q : T^*_qM \lo \ann{\D^{\text{inst}}}_q \, , \quad q
\in N \, .
\]

\subsection{Focusing points}\label{se:focusing-mechanical}

Since the mechanical system is subject to an affine distribution of
constraints, Proposition~\ref{prop:focusing-affine} implies that for a
given $u \in T^*M$, the focusing locus is $F_{(H,C)}(u) = \{ \P(u) +
\Q(\Upsilon) \}$.  Regarding the instantaneous focusing points, according
to Proposition~\ref{prop:instantaneous-focusing} one has that
$F_{(H,C^{\text{inst}},N)}(u) = \{ \P^{\text{inst}}(u) +
\Q^{\text{inst}}(\Upsilon_{C^{\text{inst}}}) \}$.

\subsection{Constrained
  characteristics}\label{se:constrained-characteristic-mechanical}

Here we give an explicit description of the characteristic curves. Let
$N$ be the critical hypersurface, and assume that (locally)
$N=f^{-1}(0)$, with $f \in C^{\infty}(M)$ verifying that $d_qf \neq 0$
for all $q \in N$. Consider the covector field $\P (df)$ along $N$
defined as $q\mapsto\P (df)_q = \P_q(d_qf), \, q\in N$.  The
transversality assumption between $C$ and $N$ implies that $\P (df)_q
\not = 0$, for all $q\in N$.  Clearly $\P (df) \in C^0$. In addition,
for $v \in \D \cap TN$,
\[
\P(df) (v) = (df - \Q (df))(v) = df (v) = 0 \, ,
\]
and one can conclude that $C^0 \cap \ann{\D \cap TN} = \spn \{ \P
(df )\}$. Therefore, we have the following result.

\begin{lemma}\label{le:characteristic-mechanical}
  The constrained characteristic of a mechanical system $(M,H,C|N)$
  passing through $y \in C \cap T^*_NM$ is given by $\zeta_y = y +
  \spn \{ \P (d_{\pi(y)}f) \} \subset C \cap T^*_NM$.
\end{lemma}

Note that in the absence of constraints one recovers the standard
non-constrained characteristic $\zeta_y = y + \spn \{d_{\pi(y)}f
\}$ passing through $y$.

\subsection{Decisive points: elastic
  mode}\label{se:mechanical-decisive-elastic}

\subsubsection{Change of
  constraints}\label{se:mechanical-decisive-elastic-change}

Let $C_{\pm} \subset T^*M$ be two affine constraint submanifolds.
Denote by $\P_\pm$ and $\Q_\pm$ the projectors corresponding to
$C_\pm$ and the co-metric $\G$.  Let $y \in C_\varepsilon \cap T^*_NM$
be a $\varepsilon$-out point, $\varepsilon \in \{+,-\}$.  Then,
according to Definition~\ref{dfn:decisive-change-constraints}, an
$y$-admissible sequence, $(y_i,\varepsilon_i)$, $i=0,1,\dots,k$, is
necessarily of the form $y_{i+1} = \P_{\varepsilon_{i+1}}(y_i) +
\Q_{\varepsilon_{i+1}} (\Upsilon_{C_{\varepsilon_{i+1}}})$. If
instantaneous constraints are present, then one has to use the
projectors $\P^{\text{inst}}_{\varepsilon}$ and
$\Q^{\text{inst}}_{\varepsilon}$ instead of $\P_{\varepsilon}$ and
$\Q_{\varepsilon}$, respectively.

\begin{remark}
  {\rm Mechanical systems subject to generalized constraints are also
    treated in~\cite{CoLeMaMa} in a somehow different context. The
    approach taken there makes use of generalized (i.e. non-constant
    rank) codistributions defining the nonholonomic constraints and a
    generalized version of Newton's second law.  Under appropriate
    regularity conditions, it can be seen that the `post-impact' point
    in~\cite{CoLeMaMa} is a decisive point of the Hamiltonian system
    according to Definition~\ref{dfn:decisive-change-constraints}.}
\end{remark}

\subsubsection{Discontinuous Hamiltonian
  systems}\label{se:mechanical-decisive-elastic-Hamiltonian}

Let $C_{\pm} \subset T^*M$ be two affine constraint submanifolds.  Let
$g_\pm$ be a Riemannian metric on $M_\pm$ and $V_\pm \in C^\infty
(M_\pm)$ such that
\begin{equation}\label{eq:various-equalities}
  H_\pm(q,p) = \hat{T}_\pm (q,p) + V_\pm(q) \, , \quad \hat{T}_\pm(q,p)
  = \frac{1}{2} \G_\pm (p,p) \, .
\end{equation}
For simplicity, we only treat the case $V_{\pm} = V|_{M_{\pm}}$, $V\in
C^{\infty}(M)$. We denote by $\P_\pm$ and $\Q_\pm$ the projectors
corresponding to $C_\pm$ and the co-metric $\G_{\pm}$.  Additionally,
let $C^{\text{inst}}_{\pm} \subset T_N^*M$ be affine constraint
submanifolds corresponding to some instantaneous constraints imposed
along $N$.  Denote by $\P_{\pm}^{\text{inst}}$ and
$\Q_{\pm}^{\text{inst}}$ the projectors corresponding
to~$C^{\text{inst}}_{\pm}$ and the co-metric~$\G_{\pm}$.

Let $y \in C_\varepsilon \cap T_N^*M$ be an $\varepsilon$-out
point. Following Definition~\ref{dfn:decisive-dis-elastic}, we
first describe the reflective and refractive steps with initial
point $(y,\varepsilon)$. According to
Lemma~\ref{le:characteristic-mechanical}, we have to look for
points of the form
\[
x = y + c \P_\varepsilon (d_qf) \, , \quad q = \pi(y) \, ,
\]
for some $c$, which in addition belong to the same $H$-energy
level as~$y$.

\paragraph{Reflective step}
Concerning the 1-st move, note that $y+c\P_\varepsilon(d_qf)$ and $y$
belong to $T_q^*M$. Then, the equality
$H_\varepsilon(y+c\P_\varepsilon(d_qf)) = H_\varepsilon(y)$ implies
that $\hat{T}_\varepsilon(y+c\P_\varepsilon(d_qf)) =
\hat{T}_\varepsilon(y)$.  Now,
\begin{multline*}
  \hat{T}_\varepsilon(y+c\P_\varepsilon(d_qf)) =\\
  \hat{T}_\varepsilon(y) + c \, \G_\varepsilon(y,\P_\varepsilon(d_qf))
  + \frac{c^2}{2}
  \G_\varepsilon(\P_\varepsilon(d_qf),\P_\varepsilon(d_qf)) \, ,
\end{multline*}
and, therefore, we have
\[
c \left( \G_\varepsilon (y,\P_\varepsilon(d_qf)) + \frac{c}{2}
  \G_\varepsilon(\P_\varepsilon(d_qf),\P_\varepsilon(d_qf)) \right) =
0 \, ,
\]
with solutions
\begin{equation}\label{eq:solutions-constrained-reflective}
  c_{\varepsilon,1} = 0 \, , \quad c_{\varepsilon,2} = - \frac{2 \,
  \G_\varepsilon(y,\P_\varepsilon(d_qf))}{\G_\varepsilon(
  \P_\varepsilon(d_qf), \P_\varepsilon(d_qf))} \, .
\end{equation}

An important property of these points is contained in the following
lemma.
\begin{lemma}\label{le:general}
  Let $y \in C_\varepsilon \cap T^*_NM$ and $c_{\varepsilon,2}$ be the
  constant given by~\eqref{eq:solutions-constrained-reflective}.
  Then,
  \begin{align*}
    \G_\varepsilon (y,d_qf) &= \G_\varepsilon (\P_\varepsilon (y) +
    \Q_\varepsilon ({\Upsilon}_q),
    d_qf) \, , \\
    \G_\varepsilon (y + c_{\varepsilon,2}
    \P_\varepsilon(d_qf),d_qf) &= \G_\varepsilon (-\P_\varepsilon (y) +
    \Q_\varepsilon ({\Upsilon}_q), d_qf) \, .
  \end{align*}
\end{lemma}

\begin{proof}
  The first statement follows by noting that if $y \in C_q$, then $y =
  \P_\varepsilon (y) + \Q_\varepsilon ( {\Upsilon}_q) $. For the
  second one, notice that
  \begin{multline*}
    \G_\varepsilon (y + c_{\varepsilon,2}
    \P_\varepsilon(d_qf),d_qf) = \G_\varepsilon (y ,d_qf) +
    c_{\varepsilon,2}
    \G_\varepsilon (\P_\varepsilon(d_qf),\P_\varepsilon(d_qf)) \\
    = \G_\varepsilon (y ,d_qf- 2 \P_\varepsilon(d_qf)) = -
    \G_\varepsilon (\P_\varepsilon(y) ,
    \P_\varepsilon(d_qf)) + \G_\varepsilon (\Q_\varepsilon(y) , d_qf) \\
    = \G_\varepsilon (-\P_\varepsilon (y) + \Q_\varepsilon (
    {\Upsilon}_q), d_qf) \, ,
  \end{multline*}
  which gives the desired result.  \qed
\end{proof}

The 2-nd move simply consists of determining the focusing points for
points~\eqref{eq:solutions-constrained-reflective} with respect to
$C^{\text{inst}}_{{\varepsilon}}$ or, if ${\varepsilon}$-instantaneous
constraints are absent, with respect to $C_{{\varepsilon}}$. This is
done in terms of the corresponding projectors, exactly as explained in
Section~\ref{se:focusing-mechanical} above.

\paragraph{Refractive step}
Concerning the 1-st move, the equality
$H_{\bar{\varepsilon}}(y+c\P_{\varepsilon}(d_qf)) = H_\varepsilon(y)$
implies $\hat{T}_{\bar{\varepsilon}}(y+c\P_{\varepsilon}(d_qf)) =
\hat{T}_{\varepsilon}(y)$.  Now,
\[
\hat{T}_{\bar{\varepsilon}}(y + c\P_{\varepsilon}(d_qf)) =
\hat{T}_{\bar{\varepsilon}}(y) + c \,
\G_{\bar{\varepsilon}}(y,\P_{\varepsilon}(d_qf)) + \frac{c^2}{2}
\G_{\bar{\varepsilon}}(\P_{\varepsilon}(d_qf),\P_{\varepsilon}(d_qf))
\, .
\]
Therefore, one has
\[
c \left( \G_{\bar{\varepsilon}} (y,\P_{\varepsilon}(d_qf)) +
  \frac{c}{2}
  \G_{\bar{\varepsilon}}(\P_{\varepsilon}(d_qf),\P_{\varepsilon}(d_qf))
\right) + \hat{T}_{\bar{\varepsilon}}(y) -\hat{T}_{\varepsilon}(y) = 0
\, ,
\]
with solutions $i=1,2$,
\begin{multline}\label{eq:solutions-constrained-refractive}
  c_{\bar{\varepsilon},i} =
  \frac{1}{\G_{\bar{\varepsilon}}(\P_{\varepsilon}(d_qf)),
    \P_{\varepsilon}(d_qf))} \Big(-
  \G_{\bar{\varepsilon}}(y,\P_{\varepsilon}(d_qf)) \pm \\
  \sqrt{\G_{\bar{\varepsilon}}(y,\P_{\varepsilon}(d_qf))^2-2
    \G_{\bar{\varepsilon}}(\P_{\varepsilon}(d_qf),\P_{\varepsilon}(d_qf))
    (\hat{T}_{\bar{\varepsilon}}(y) -\hat{T}_{\varepsilon}(y))} \Big)
  \, .
\end{multline}


As before, the 2-nd move simply consists of computing the focusing
points for the solutions~\eqref{eq:solutions-constrained-refractive}
with regards to $C^{\text{inst}}_{\bar{\varepsilon}}$ or, if
$\bar{\varepsilon}$-instantaneous constraints are absent,
$C_{\bar{\varepsilon}}$. This is done in terms of the corresponding
projectors according to Section~\ref{se:focusing-mechanical}.


\subsubsection*{Discontinuous Hamiltonian systems with smooth
  constraints}

In this situation, there is a single constraint submanifold $C$, and a
discontinuous Hamiltonian $H_\pm$ on $T^*M$.  Denote by $\P_\pm$ and
$\Q_\pm$ the projectors corresponding to $C$ and the
co-metrics~$\G_\pm$, respectively. According to
Proposition~\ref{prop:decisive-change-Hamiltonian}, the decisive
points for a given $\varepsilon$-out point $y\in C \cap T_N^*M$ are
simply the in and trapping points belonging to the intersection of the
constrained characteristic $\zeta_y$ passing through $y$ with the set
$\setdef{z\in C}{H_{\pm}(z)=H_{\varepsilon}(y)}$.  Therefore, as
candidate $\varepsilon$-decisive points we have the solution
corresponding to $c_{\varepsilon,2}$
in~\eqref{eq:solutions-constrained-reflective}, and as candidate
$\bar{\varepsilon}$-decisive points we have the solutions
corresponding to $c_{\bar{\varepsilon},i}$, $i=1,2$,
in~\eqref{eq:solutions-constrained-refractive}.


\begin{proposition}\label{prop:duality}
  Let $y \in C_\varepsilon \cap T_N^*M$ be a $\varepsilon$-out point.
  If the constraints are linear, $C=C^0$, then the solution
  corresponding to $c_{\varepsilon,2}$
  in~\eqref{eq:solutions-constrained-reflective} is a
  $\varepsilon$-decisive point for $y$.
\end{proposition}

\begin{proof}
  The basic observation is the second order character of the dynamics,
  both in the presence and in the absence of nonholonomic constraints.
  This implies that for any $y \in T^*M$ and any distribution of
  affine constraints $C$, we have $X_{H}(f)(y) = X_{H,C}(f)(y)$, since
  $f$ is only a function of the configurations.  Note that if $H$ is
  of mechanical type, then $X_{H} (f) (x) = \G(x,d_qf)$, for any $x
  \in T_q^*M$.  Now, from Lemma~\ref{le:general}, taking
  $H=H_\varepsilon$, one gets
  \[
  \G_\varepsilon(y + c_2 \P_\varepsilon(d_qf),d_qf) = -
  \G_\varepsilon(y,d_qf) \, .
  \]
  Since $y$ is a $\varepsilon$-out point, then $\G_\varepsilon
  (y,d_qf) \neq 0$. Consequently, $X_{H_+} (f) (y + c_2 \P_+(d_qf)) =
  - \G_+(d_qf,y)$ has the opposite sign, and hence it is an in point.
  \qed
\end{proof}






\subsection{Decisive points: inelastic
  mode}\label{se:mechanical-decisive-inelastic}

\subsubsection{Change of
  constraints}\label{se:mechanical-decisive-inelastic-change}

Let $C_{\pm} \subset T^*M$ be two affine constraint submanifolds and
let $C^{\text{inst}}$ be a set of instantaneous affine constraints. We
denote by $\P^{\text{inst}}_\pm$ and $\Q^{\text{inst}}_\pm$ the
projectors corresponding to $C^{\text{inst}}_\pm$ and the co-metric
$\G$. If the instantaneous constraints are absent, denote by $\P_\pm$
and $\Q_\pm$ the projectors corresponding to $C^{\text{tr}}_{\pm}$ and
the co-metric $\G$.  Let $y \in C_\varepsilon \cap T_N^*M$ be a
$\varepsilon$-out point, $\varepsilon \in \{+,-\}$.  Then, according
to Definition~\ref{dfn:decisive-change-inelastic}, the unique
$y$-decisive point is $\P^{\text{inst}}_{\bar{\varepsilon}}(y) +
\Q^{\text{inst}}_{\bar{\varepsilon}}
(\Upsilon_{C_{\bar{\varepsilon}}})$ (or, if there are no instantaneous
constraints, $\P_{\bar{\varepsilon}}(y) + \Q_{\bar{\varepsilon}}
(\Upsilon_{C_{\bar{\varepsilon}}})$).

\subsubsection{Discontinuous Hamiltonian
  systems}\label{se:mechanical-decisive-inelastic-Hamiltonian}

As in Section~\ref{se:mechanical-decisive-elastic-Hamiltonian}, let
$C_{\pm} \subset T^*M$ be two affine constraint submanifolds, $g_\pm$
a Riemannian metric on $M_\pm$ and $V_\pm \in C^\infty (M_\pm)$ such
that equation~\eqref{eq:various-equalities} is verified.  For
simplicity, we only treat the case $V_{\pm} = V|_{M_{\pm}}$, $V\in
C^{\infty}(M)$. We denote by $\P_\pm$ and $\Q_\pm$ the projectors
corresponding to $C_\pm$ and the co-metric $\G_{\pm}$. Additionally,
let $C^{\text{inst}}_{\pm} \subset T_N^*M$ be affine constraint
submanifolds corresponding to some instantaneous constraints imposed
along $N$.  We denote by $\P^{\text{inst,tr}}_{\varepsilon}$ and
$\Q^{\text{inst,tr}}_{\varepsilon}$ the projectors associated with the
submanifold $C^{\text{inst,tr}}_{\varepsilon}$ and the co-metric
$\G_{\varepsilon}$. In the absence of instantaneous constraints, we
denote by $\P^{\text{tr}}_{\varepsilon}$ and
$\Q^{\text{tr}}_{\varepsilon}$ the projectors associated with the
submanifold $C^{\text{tr}}_{\varepsilon}$ and the co-metric
$\G_{\varepsilon}$. In case $N = \partial M$, we denote the latter
with the superscript ``$\partial$'' instead of ``tr''.

Let $y \in C_\varepsilon\cap T_N^*M$ be an $\varepsilon$-out point.
The points associated with $y$ resulting from the 1-st moves in a
reflected or a refracted falling step are given, respectively, by
equations~\eqref{eq:solutions-constrained-reflective}
and~\eqref{eq:solutions-constrained-refractive}.  As before, the 2-nd
move simply consists of computing the focusing points for these
solutions with respect to $C^{\text{inst,tr}}_{\varepsilon}$ for a
reflected falling step (respectively,
$C^{\text{inst,tr}}_{\bar{\varepsilon}}$ for a refracted falling step)
or, if the instantaneous constraints are absent,
$C^{\text{tr}}_\varepsilon$ (respectively,
$C^{\text{tr}}_{\bar{\varepsilon}}$).  This is done in terms of the
corresponding projectors according to
Section~\ref{se:focusing-mechanical}. According to
Definition~\ref{dfn:decisive-dis-inelastic}, this gives all the
$y$-decisive points.

\begin{proposition}\label{prop:decisive-dis-inelastic}
  Let $y \in C_\varepsilon\cap T_N^*M$ be an $\varepsilon$-out point
  and assume that the constraints are linear. For $N = \partial M$,
  the unique $y$-reflected falling point is given by
  $\P^{\text{inst},\partial}_{\varepsilon} (y)$ (or, in the absence of
  $\varepsilon$-instantaneous nonholonomic constraints,
  $\P^{\partial}_{\varepsilon} (y)$).
\end{proposition}

\begin{proof}
  From the previous discussion, we know that the points in the
  constrained characteristic passing through $y$ with the same
  $H_{\varepsilon}$-energy level are $y$ itself and $y +
  c_{\varepsilon,2} \P_\varepsilon (d_qf)$, $q = \pi(y)$ (cf.
  equation~\eqref{eq:solutions-constrained-reflective}). Now, note
  that $d_qf$ belongs to the $\G_\varepsilon$-orthogonal complement of
  $\alpha_{H_\varepsilon}^{-1} (T(\partial M))$, i.e.
  \begin{align*}
    \G_\varepsilon (d_qf, \beta) = d_qf (\alpha_{H_\varepsilon} (\beta)) =
    0 \, , \quad \beta \in \alpha_{H_\varepsilon}^{-1} (T(\partial M))
    \, .
  \end{align*}
  Using the equality $C^{\partial}_{\varepsilon} = C_{\varepsilon}
  \cap \alpha_{H_\varepsilon}^{-1} (T(\partial M))$, we have that $d_qf
  \in \alpha_{H_\varepsilon}^{-1} (T(\partial M))^{\perp_\varepsilon}$
  implies $\P^{\text{inst},\partial}_{\varepsilon} (\P_\varepsilon
  (d_qf)) = 0$ and $\P^{\partial}_{\varepsilon} (\P_\varepsilon (d_qf)) =
  0$.  The result is then a consequence of
  Definition~\ref{dfn:decisive-dis-inelastic}.  \qed
\end{proof}

\subsection{Energy behavior}

In this section, we discuss the consequences regarding the energy
behavior of the system that result from the application of the
Transition Principle.

\begin{lemma}\label{le:simple}
  Given $y \in T^*_NM$, let $x=\P(y) + \Q({\Upsilon}_q)$, $q =
  \pi(y),$ be the associated $y$-focusing point with respect to a
  submanifold $C \subset T^*M$. Then
  \[
  \hat{T}(x) \le \hat{T}(y) +
  \hat{T}(\Q({\Upsilon}_q),\Q({\Upsilon}_q)) \, ,
  \]
  and the equality holds if and only if $y$ belongs to $C^0$.
\end{lemma}

\begin{proof}
  Note that
  \begin{multline*}
    \G (\P(y)+\Q({\Upsilon}_q),\P(y)+\Q({\Upsilon}_q)) = \G
    (\P(y),\P(y)) + \G (\Q({\Upsilon}_q),\Q({\Upsilon}_q))) \\
    \le \G (\P(y),\P(y)) + \G (\Q(y),\Q(y)) + \G
    (\Q({\Upsilon}_q),\Q({\Upsilon}_q))\\
    = \G (y,y) + \G (\Q({\Upsilon}_q),\Q({\Upsilon}_q)) \, ,
  \end{multline*}
  where we have used that $\G$ is positive-definite, and the fact that
  $C^0$ and $\ann{\D}$ are orthogonal spaces with respect to the
  co-metric $\G$.  If the equality holds, then $\G (\Q(y),\Q(y))= 0$,
  which is equivalent to $y \in C^0$.  \qed
\end{proof}

As a consequence of this simple lemma we can conclude that in the case
of linear constraints the Transition Principle always implies a loss
of energy.  This is a suitable generalization to constrained systems
of the classical Carnot theorem for systems subject to impulsive
forces~\cite{Ro}.

\begin{theorem}[Carnot's theorem for generalized linear
  constraints] Suppose that the Hamiltonian system is subject to
  nonholonomic constraints given by a linear distribution. Then the
  Transition Principle implies always a loss of energy as the result
  of an ``impact''.
\end{theorem}

\begin{proof}
  Under linear constraints, note that $\Upsilon=0$. From
  Lemma~\ref{le:simple}, we get $\hat{T}(x) \le \hat{T}(y)$ with
  $F_{(H,C)}(y) = \{ x \}$. The result now follows from the
  formulation of the Transition Principle and the definitions of
  decisive points in Section~\ref{se:decisive} (cf.
  Definitions~\ref{dfn:decisive-change-constraints}-\ref{dfn:decisive-dis-inelastic}).
  \qed
\end{proof}

Under linear constraints, the trajectory of the system maintains the
same energy level after the application of the Transition Principle in
the following cases:
\begin{enumerate}
\item when the decisive points are determined according to
  Definitions~\ref{dfn:decisive-change-constraints}
  and~\ref{dfn:decisive-change-inelastic} and the impact point $y \in
  T^*_NM$ belongs to $C_+ \cap C_-$; and
\item when the constraints are smooth, and therefore the decisive
  points are determined according to
  Proposition~\ref{prop:decisive-change-Hamiltonian}.
\end{enumerate}
If the decisive points are determined according to
Definitions~\ref{dfn:decisive-dis-elastic}
and~\ref{dfn:decisive-dis-inelastic}, then nothing can be said in
general. The refractive steps will typically imply an energy decrease.

\begin{remark}
  {\rm This type of energy arguments also allows to discard as follows
    the possibility of chattering when computing the $y$-decisive
    points if the constraints change (see
    Definition~\ref{dfn:decisive-change-constraints} and
    Remark~\ref{rem:chattering} above).  Let $N=\setdef{y \in
      T^*M}{f(y)=0}$.  Assume there is an infinite $y$-admissible
    sequence $(y_i,\varepsilon_i)$, $i=0,\dots,\infty$.  For each $i$,
    we have that $y_i \neq y_{i+1}$, since otherwise
    \[
    X^l_{H,C_{\varepsilon_{i+1}}} (f) (y_{i+1}) = X^l_{H} (f)
    (y_{i+1}) = X^l_{H,C_{\varepsilon_{i}}} (f) (y_{i}) \, , \quad
    \text{for all $l$} \, ,
    \]
    which together with the fact that $y_i$ is a $\varepsilon_i$-out
    point, implies that $y_{i+1}$ is a $\varepsilon_{i+1}$-in point.
    The latter contradicts the definition of admissible sequence. As a
    consequence of Lemma~\ref{le:simple}, we then have
    \begin{align*}
      \hat{T}(y_0) > \hat{T}(y_1) > \hat{T}(y_2) > \dots >
      \hat{T}(y_i) \ge 0 \, .
    \end{align*}
    The limit of this sequence is zero, which implies that the
    $y$-decisive point corresponding to such a sequence would be $0$,
    that is, the trajectory would get `stuck' when reaching~$N$.  }
\end{remark}

\subsection{Integrable constraints}

The integrability of the constraints simplifies the application of the
Transition Principle. Consider, for instance, the situation when the
mechanical system is unconstrained on $M_-$ and is subject to some
generalized linear constraints $C = C^0$ on $M_+$ that turn out to be
holonomic, i.e., $\alpha_H(C^0) = \D$ is integrable. Denote by
$\{S_\alpha \}$, $\alpha$ being an $m$-dimensional parameter, the
foliation of $M_+$ induced by $\D$.  Locally this foliation is
described by $m$ functions $f_i\in C_\infty(M)$ such that
\[
q \in S_\alpha \Longleftrightarrow f_i(q)=\alpha_i \, , \; 1 \le i \le
m \, .
\]
A similar situation has been treated in~\cite{PuVi2} in the
context of totally inelastic collisions (note, however, that
in~\cite{PuVi2} the integrable distribution is defined only on
$N$, whereas here $\D$ is defined on $M_+$).  The integrable
constraints imposed by $\D$ can be interpreted as an abrupt
reduction of the phase space of the mechanical system.

By definition, one has that $\ann{\D} =\spn \{ df_1, \dots, df_m \}$.
The matrix $J$ in~\eqref{eq:matrices} is then given by $J=(\partial
f_i / \partial q^a)$ and the projector $\P$ is $\P (x) = (1-J^t
\mathcal{B}^{-1} J G^{-1})x$. Let $y \in C_- \cap T^*_NM$ be the
impact state of a trajectory $q(t)$ coming from $M_-$. From the
discussion in Section~\ref{se:mechanical-decisive-elastic}, we obtain
that the unique focusing point associated to $y$ is $ x= (1-J^t
\mathcal{B}^{-1} J G^{-1})y$. The trajectory will continue its motion
in $M_+$, $M_-$ or $N$ depending on the in/out/trapping character of
the focusing point $x$.  If it evolves in $M_+$ (more precisely, in
$S_\alpha \subset M_+$ with $\alpha$ such that $x\in S_\alpha$), we
call it the `refraction' of the original trajectory. If it evolves in
$M_-$, we call it the `reflection' of the original trajectory.


\section{Examples}\label{se:examples}

In this section we consider four examples to illustrate the theory
exposed above.  They all present the example of a rolling sphere
considered in various constrained situations. The first one is taken
from~\cite{CoLeMaMa} and is treated here in order to provide a further
comparison with previous approaches. The second one combines the
presence of smooth nonholonomic constraints with discontinuous
Hamiltonians and instantaneous constraints acting on the system along
a hypersurface.  The third one consists of a ball rolling on a
rotating surface whose angular velocity is suddenly changed to a
different value, and this is modeled via a discontinuous affine
distribution of constraints.  Finally, the fourth one presents a
two-wheeled system with a rod of variable length and illustrates the
application of the Transition Principle in both the elastic and the
inelastic modes.

\subsection{A rolling sphere}\label{se:rolling-sphere}

\begin{figure}[htb!]
  \vspace*{-.5cm}
  \begin{center}
    \includegraphics[width=8cm,height=4.5cm]{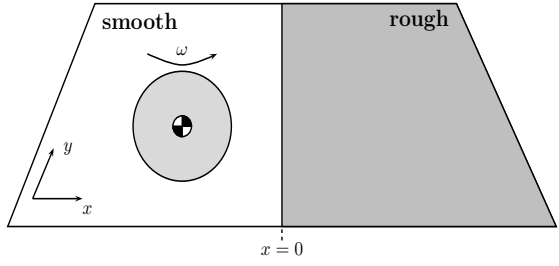}
  \end{center}
  \vspace*{-.5cm}
  \caption{The rolling sphere on a `special' surface.}\label{fig:bola}
\end{figure}

Consider a homogeneous sphere rolling on a plane.  Assume it has unit
mass ($m=1$) and let $k^2$ be its inertia about any axis. Let $(x,y)$
denote the position of the center of the sphere and let $(\varphi,
\theta, \psi)$ denote the Eulerian angles.  The configuration space is
therefore $Q=\real^2 \times \hbox{SO}(3)$.  Assume that the plane is
smooth if $x<0$ and absolutely rough if $x>0$ (see
Figure~\ref{fig:bola}). On the smooth half-plane, the motion of the
sphere is assumed free, that is, the sphere can slip. On the rough
half-plane, the sphere should roll without slipping due to the
constraints imposed by the roughness.  We are interested in
determining the eventual sudden changes in the trajectories of the
sphere when it reaches the line separating the smooth and the rough
half-planes.

The kinetic energy of the sphere is
\begin{equation}\label{eq:Lagrangian}
  T = \frac{1}{2}
  \left(\dot{x}^2+\dot{y}^2+k^2(\omega_x^2+\omega_y^2+
    \omega_z^2)\right) \, ,
\end{equation}
where $\omega_{x}$, $\omega_{y}$ and $\omega_{z}$ are the angular velocities
with respect to the inertial frame, given by
\begin{align*}
  \omega_{x} = \dot{\theta} \cos \psi + \dot{\varphi} \sin \theta \sin
  \psi \, , \quad \omega_{y} = \dot{\theta} \sin \psi - \dot{\varphi}
  \sin \theta \cos \psi \, , \quad \omega_{z} = \dot{\varphi} \cos
  \theta + \dot{\psi} \, .
\end{align*}

The condition of rolling without sliding of the sphere when $x>0$
implies that the point of contact of the sphere and the plane has zero
velocity
\begin{eqnarray*}
  \phi^1 = \dot{x}-r\omega_y=0 \, , \quad \phi^2 = \dot{y}+r\omega_x=0 \, ,
\end{eqnarray*}
where $r$ is the radius of the sphere.

Following the classical procedure~\cite{NeFu}, we introduce
quasi-coordinates `$q^{1}$', `$q^{2}$' and `$q^{3}$' such that
$\dot{q}^{1}$'$ = \omega_x$, `$\dot{q}^{2}$'$= \omega_y$ and
`$\dot{q}^{3}$'$=\omega_z$. The latter merely have a symbolic meaning
in the sense that in the present example, for instance, the partial
derivative operators $\partial/\partial q^i$ should be interpreted as
linear combinations of the partial derivatives with respect to Euler's
angles.  Also to the differential forms $dq^i$ one should attach the
appropriate meaning, i.e.  they do not represent exact differentials
but, instead, we should read them as $dq^1 = \cos\psi\,d\theta
+\sin\theta\sin\psi\,d\varphi$, etc.

The singular hypersurface $N$ is defined by $N = \{ x=0 \}$. In this
case, the constraints are linear and the nonholonomic distribution
$\alpha_H(C)=\D$ on $M_+$ is given by
\[
\D_{(x,y,q^1,q^2,q^3)} = \hbox{span} \left\{ r \pder{}{x} +
  \pder{}{q^2} , - r \pder{}{y} + \pder{}{q^1}, \pder{}{q^3} \right\}
\, .
\]
Here we are dealing with a single distribution which constrains the
motion on~$M_+$.

In the following we compute the decisive points for this example.  Let
$\lambda \in C_- \cap T^*M$ be a minus-out point. A direct computation
shows that the expression of the projector $\P:T^*M \rightarrow C$ in
local coordinates is
\begin{align}\label{eq:projector}
  \P = \left(
    \begin{array}{ccccc}
      \frac{r^{2}}{r^{2}+k^{2}}&0&0&\frac{r}{r^{2}+k^{2}}&0\\
      0&\frac{r^{2}}{r^{2}+k^{2}}&-\frac{r}{r^{2}+k^{2}}&0&0\\
      0&\frac{-rk^{2}}{r^{2}+k^{2}}&\frac{k^{2}}{r^{2}+k^{2}}&0&0\\
      \frac{rk^{2}}{r^{2}+k^{2}}&0&0&\frac{k^{2}}{r^{2}+k^{2}}&0\\
      0&0&0&0&1
    \end{array}
  \right) \, .
\end{align}
Therefore, the single focusing point for $\lambda \in T_N^*M$ is given
by $x = \P (\lambda ) \in C \cap T_N^*M$. If we denote $\lambda
=(x_0,y_0,q^1_0,q^2_0,q^3_0,(p_x)_0,(p_y)_0,(p_1)_0,(p_2)_0,(p_3)_0)$
and $x = (x,y,q^1,q^2,q^3,p_x,p_y,p_1,p_2,p_3)$, we get
\[
\setlength{\arraycolsep}{2pt}
\begin{array}{rcl}
  p_{x} & = & \displaystyle{\frac{r^{2} (p_{x})_{0} +
      r (p_{2})_{0}}{r^{2}+k^{2}}} \; ,\\
  p_{y} & = & \displaystyle{\frac{r^{2} (p_{y})_{0} -
      r (p_{1})_{0}}{r^{2}+k^{2}}} \; ,
\end{array}
\quad
\begin{array}{rcl}
  p_{1} & = & \displaystyle{\frac{-rk^{2} (p_{y})_{0} +
      k^{2} (p_{1})_{0}}{r^{2}+k^{2}}} \; ,\\
  p_{2} & = & \displaystyle{\frac{rk^{2} (p_{x})_{0} +
      k^{2} (p_{2})_{0}}{r^{2}+k^{2}}} \; ,\\
  p_{3} & = & (p_{3})_{0} \; .
\end{array}
\]
Note also that the focusing point with respect to $C_- = T^*M$
associated with $x$ is $x$ itself. Therefore, if $x$ is a plus-out
point, the only admissible sequence for $\lambda$ is $\{(\lambda,-),
(x,+), (x,-)\}$. If $x$ is either a plus-in or a plus-trapping point,
then the only admissible sequence for $\lambda$ is $\{(\lambda,-),
(x,+)\}$.  The set of plus-trapping points for the dynamics
$X_{H,C_+}$ is $\partial (T^*M)^n = \setdef{\mu \in T^*M}{x=0 \, , \;
  p_x = 0}$.  Consequently, the trajectory is refracted, i.e., the
sphere follows its motion on $M_+$ under the dynamics $X_{H,C_+}$
(rolling without slipping) if $p_x \ge 0$. Otherwise (i.e., if $p_x <
0$), the trajectory is reflected by the ``roughness'' and continues in
$M_-$ under the dynamics $X_H$ starting from $x$.

\subsection{A rolling sphere hitting a wall}

\begin{figure}[htb!]
  \vspace*{-.5cm}
  \begin{center}
    \includegraphics[width=8cm,height=6.5cm]{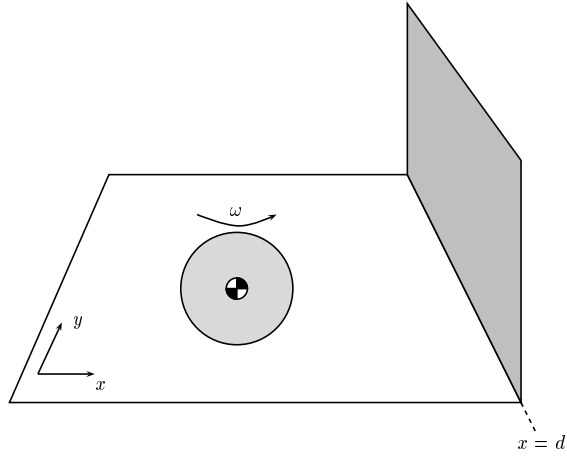}
  \end{center}
  \vspace*{-.5cm}
  \caption{A rolling sphere that eventually hits a
    wall.}\label{fig:ball-hitting}
\end{figure}

This is a classical example~\cite{CeLaRe,IbLeLaMaMaPi,NeFu} that we
treat here for the sake of completeness.  Consider again a homogeneous
sphere of radius $r$ and unit mass. Assume that the sphere rolls
without sliding on a horizontal table, and that at a certain instant
of time it hits a wall determined by the plane $x=d >0$ (cf.
Figure~\ref{fig:ball-hitting}).  When this happens, the following
constraint is instantaneously imposed on the system,
\begin{align*}
  \psi = \dot{y} - r \omega_z = 0 \, . 
\end{align*}
Therefore, we are in the situation explained in
Remark~\ref{rem:decisive-boundary}.
The configuration space of the system is $M=M_+=\{x<d\}$, with the
boundary $N=\{x=d\}$, and the linear constraint submanifold
$C=C_+\subset T^*M$ is given by $\alpha_H (C) = \D$,
\begin{align*}
  \D & = \text{span} \left\{ r \pder{}{x} + \pder{}{q^2} , - r
    \pder{}{y} + \pder{}{q^1}, \pder{}{q^3} \right\} \, .
\end{align*}
The expression for the projector $\P:T^*M \rightarrow C$ is given by
equation~\eqref{eq:projector}.  The submanifold giving the
instantaneous constraints along $N$ is
\begin{align*}
  C^{\text{inst}} = \setdef{\lambda \in C}{\psi(\alpha_H(\lambda)) =
    0} \, .
\end{align*}
The projector $\P^{\text{inst}} = \P_{+}^{\text{inst}} : T^*M
\rightarrow C^{\text{inst}}$ is
\begin{multline}\label{eq:projector-+}
  \P^{\text{inst}} (\lambda)= \frac{r \lambda_x + \lambda_2}{r^2 +
    k^2} \left( r d{x} +k^2 d{q^2} \right) \\
  + \frac{-r \lambda_y + \lambda_1-\lambda_3}{r^2 + 2 k^2} \left( -r
    d{y} +k^2 d{q^1} - k^2d{q^3} \right) \, .
\end{multline}

Let $\lambda = (x_0,y_0,q^1_0,q^2_0,q^3_0,
(p_x)_0,(p_y)_0,(p_1)_0,(p_2)_0,(p_3)_0) \in C_+ \cap T_N^*M \subset
T^*M$ be a plus-out point, i.e., $ \G (\lambda,dx) <0$. We first
consider an elastic impact. Since $H_- = \infty$, we only compute the
outcome of a reflective step.  According
to~\eqref{eq:solutions-constrained-reflective}, the points in the
constrained characteristic passing through $\lambda$ within the same
$H_+$-energy level are
\begin{align*}
  \lambda \quad \text{and} \quad \lambda + c_{+,2} \P (dx) \, , \quad
  \text{with} \; c_{+,2} = -2 \, \frac{r^2+k^2}{r^2} (p_x)_0 \, .
\end{align*}
The associated focusing points are given by
\begin{align}\label{eq:both}
  \P^{\text{inst}} (\lambda) \quad \text{and} \quad
  \P^{\text{inst}}(\lambda) + c_{+,2} \P^{\text{inst}}(\P (dx)) \, .
\end{align}
Note that $\P^{\text{inst}}(\P (dx)) = \P^{\text{inst}}(dx) = \P(dx)$,
and therefore the points in~\eqref{eq:both} belong to the same
constrained characteristic and to the same $H_+$-energy level.
Denoting the coordinates of the point $\P^{\text{inst}} (\lambda)$ by
$(x,y,q^1,q^2,q^3, p_x,p_y,p_1,p_2,p_3)$, we get
\[
\setlength{\arraycolsep}{2pt}
\begin{array}{rcl}
  p_{x} & = & (p_x)_0 \\
  p_{y} & = & \displaystyle{\frac{ (r^2 +k^2) (p_{y})_{0} +
      r (p_{3})_{0}}{r^{2}+2 k^{2}}} \; ,
\end{array}
\quad
\begin{array}{rcl}
  p_{1} & = & -k^2 \displaystyle{\frac{ (r^2 +k^2) (p_{y})_{0} +
      r (p_{3})_{0}}{r(r^{2}+2 k^{2})}} \\
  p_{2} & = & k^2 \displaystyle{\frac{(p_x)_0}{r}} \\
  p_{3} & = & k^2 \displaystyle{\frac{ (r^2 +k^2) (p_{y})_{0} +
      r (p_{3})_{0}}{r(r^{2}+2 k^{2})}} \; ,
\end{array}
\]

Now, notice that $\P^{\text{inst}}(\lambda)$ is a plus-out point,
because $\G (\P^{\text{inst}}(\lambda),dx) = \G (\lambda,dx) <0$.
Therefore, following Proposition~\ref{prop:duality}, we conclude that
the sequence $\{(\lambda,+),\P^{\text{inst}}(\lambda) + c_{+,2}
\P^{\text{inst}} (dx),+) \}$ is $\lambda$-admissible, and
$\P^{\text{inst}}(\lambda) + c_{+,2} \P^{\text{inst}} (dx)$ is a
decisive point. The other possible $\lambda$-admissible sequence
corresponds to
\[
\{(\lambda,+), \P^{\text{inst}}(\lambda), \P^{\text{inst}}(\lambda) +
c_{+,2} \P^{\text{inst}} (dx),+) \} \, ,
\]
but renders the same decisive point.

In the case of an inelastic impact,
Proposition~\ref{prop:decisive-dis-inelastic} yields that the
unique $\lambda$-decisive point is
$\P^{\text{inst},\partial}(\lambda)$. After the impact, the ball
continues its motion along the wall under the dynamics specified
by the vector field $X^{\text{inst},\partial}_+$.

\subsection{A rolling sphere  on a rotating
  table}\label{se:homogeneous-ball}

\begin{figure}[htb!]
  \vspace*{-.5cm}
  \begin{center}
    \includegraphics[width=8cm,height=4cm]{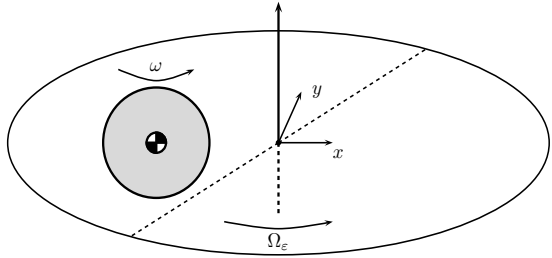}
  \end{center}
  \vspace*{-.5cm}
  \caption{A rolling sphere on a rotating
    table.}\label{fig:ball-rotating-table}
\end{figure}

Consider again a homogeneous sphere of radius $r$ and unit mass.
Assume that the sphere rolls without sliding on a horizontal table
which is rotating with certain constant angular velocity about a
vertical axis through one of its points (see
Figure~\ref{fig:ball-rotating-table}).  Let $\Omega_-$ and $\Omega_+$
be two angular velocities.  Here we consider the following situation:
each time the sphere reaches the hypersurface $x=y$, an impulsive
force is exerted on the table to put it spinning with a different
angular velocity. That is, if the angular velocity of the table was
$\Omega_-$, the force applied on its rotation axis changes it to
$\Omega_+$ and vice versa.  We assume that $\Omega_- < \Omega_+$. This
can be modeled as thinking of a system which is subject to two
different affine constraint distributions. In order for this model to
be consistent, we also assume that the surface of the table is rough
enough so that sphere is rolling without slipping at all times.

The Lagrangian is again given by equation~\eqref{eq:Lagrangian}.  The
nonholonomic constraints are now affine in the velocities,
\[
\dot{x} - r\omega_y = -\Omega y, \qquad \dot{y}+r\omega_x=\Omega x \,.
\]
The constraint space $\alpha_H(C)$ is then described by
\[
\alpha_H(C) = \D + Y = \spn \left\{ r \pder{}{x} +\pder{}{q^2}, - r
  \pder{}{y} + \pder{}{q^1}, \pder{}{q^3} \right\} + Y \, ,
\]
where $Y$ is the vector field defined by
\[
Y = - \Omega y \pder{}{x} + \Omega x \pder{}{y} \, .
\]
Note that the projection of $\Upsilon = \Legendre_L (Y)$ to
$\ann{\D}$ is given by
\[
\Q (\Upsilon) = \frac{\Omega k^2}{r^2+k^2} \left( -y dx + x dy - x r
  dq^1 - y r dq^2 \right) \, .
\]
Following the discussion for the case of affine constraints, given $y
\in T^*M$, the focusing point with respect to $C_\pm$ is given by
\[
x = \P (y) + \Q (\Upsilon_{\pm}) \, ,
\]
where $\P$ is the projector in~\eqref{eq:projector}.

Assume that the sphere is rolling on the hyperplane $M_-=\{ x < y
\}$ and that the constant angular velocity of the table is
$\Omega_-$.  Consider the case when the sphere ``hits'' the
hypersurface $N=\{ x=y \}$ with the impact state
\[
\lambda =
(x_0,y_0,q^1_0,q^2_0,q^3_0,(p_x)_0,(p_y)_0,(p_1)_0,(p_2)_0,(p_3)_0)
\in C_- = C_-^o + \Upsilon_{-} \, .
\]
Denote the coordinates of the associated focusing point by
\begin{align*}
  x = \P (\lambda) + \Q (\Upsilon_{+}) =
  (x,y,q^1,q^2,q^3,p_x,p_y,p_1,p_2,p_3) \,.
\end{align*}
Then
\begin{align}\label{eq:det-character}
  \G (df, x) = p_x - p_y = (p_x)_0 - (p_y)_0 +
  \frac{k^2}{r^2+k^2}(x_0+y_0)(\Omega_- - \Omega_+) \, .
\end{align}
Given that $\lambda$ is an minus-out point, we have that $\G
(df,\lambda) = (p_x)_0 - (p_y)_0 > 0$. If $x_0=y_0 < 0$, then the
second term in~\eqref{eq:det-character} is also positive, and
$\{(\lambda,-),(x,+)\}$ is the unique admissible sequence for
$\lambda$. In this case, $x$ is the $\lambda$-decisive point.  On the
contrary, for certain values of $x_0=y_0 > 0$, it might happen that
$\G (df, x)$ is negative, i.e., that $x$ is a plus-out point. Now,
note that the focusing point associated with $x$ is $\lambda$ itself,
since
\begin{align*}
  \P (x) + \Q (\Upsilon_{-}) = \P (\P (\lambda)) + \Q (\Upsilon_{-}) =
  \P (\lambda) + \Q (\Upsilon_{-}) = \lambda \, .
\end{align*}
As a consequence, in this case there would not be any
$\lambda$-decisive point. This problem stems from the fact the
modeling of this example as a system subject to affine constraints
does not take into account that the jump in the angular velocity of
the table takes place \emph{no matter what}.  Therefore, after the
impact, we should really regard $C_+$ as the new set of affine
constraints acting on the whole configuration manifold.  With this
interpretation, $x$ would obviously be a plus-in point (and hence
decisive). In other words, the trajectory of the ball gets reflected
back by the blow of the greater velocity $\Omega_+$.

\subsection{A two-wheeled system with a rod of variable
  length}\label{se:wheels-with-telescopic-rid}

\begin{figure}[htb!]
  \vspace*{-.4cm}
  \begin{center}
    \includegraphics[width=8.5cm]{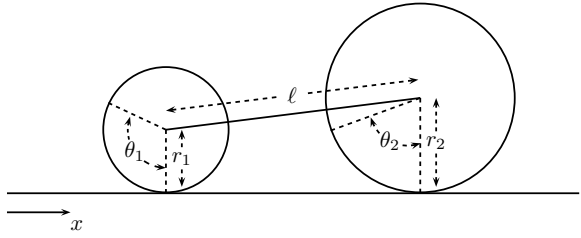}
  \end{center}
  \vspace*{-.4cm}
  \caption{A two-wheeled system with a rod of variable
    length.}\label{fig:wheels-with-telescopic-rod}
\end{figure}

Consider a system composed of two wheels of different radii, $r_1
<r_2$, connected by a massless rod of variable length $\ell$ (see
Figure~\ref{fig:wheels-with-telescopic-rod}).  For simplicity,
assume that the two-wheeled system moves along a line, and that
both the masses and the momenta of inertia of the wheels are
unitary. The wheels are subject to the standard constraints of
non-slipping. Assume that the length $\ell$ of the rod is
constrained between a minimum length $a$ and a maximum length $b$.
Here we consider the following two situations: (i) when the length
$\ell$ of the rod becomes extreme, an elastic impact occurs; (ii)
when the length $\ell$ of the rod becomes extreme, an arresting
device fixes it, and therefore an inelastic impact occurs.

The Lagrangian of the system is given by the kinetic energy of the
wheels
\begin{align*}
  L = \frac{1}{2} \left( \dot{\theta}_1^2 + \dot{\theta}_2^2 +
    \dot{x}_1^2 + \dot{x}_2^2 \right) \, .
\end{align*}
The conditions of rolling without sliding are encoded in the
constraints
\begin{align*}
  \dot{x}_1 - r_1 \dot{\theta}_1 = 0 \, , \quad \dot{x}_2 - r_2 \,
  \dot{\theta}_2 = 0 \, ,
\end{align*}
which, since we are considering the motion of the two-wheeled system
only along a line, turn out to be holonomic. The constraint on the
length of the rod is given by
\begin{align*}
  a \le \ell = \sqrt{(r_2-r_1)^2 + (x_2-x_1)^2} \le b \, .
\end{align*}
Following Remark~\ref{rem:decisive-boundary}, we set $M_- =
\emptyset$, $M_+ = M= \setdef{(x_1,x_2,\theta_1,\theta_2) \in \real^2
  \times \sphere^1\times \sphere^1}{a \le
  \ell(x_1,x_2,\theta_1,\theta_2) \le b}$, with boundary set $N =
\partial M = \setdef{(x_1,x_2,\theta_1,\theta_2) \in \real^2 \times
  \sphere^1\times \sphere^1}{\ell(x_1,x_2,\theta_1,\theta_2) = a \;
  \text{or} \; \ell(x_1,x_2,\theta_1,\theta_2) = b}$, and linear
constraint submanifold $C = C_+ \subset T^*M$ given by $\alpha_H (C) =
\D$,
\begin{align*}
  \D = \text{span} \left\{ r_1 \pder{}{x_1} + \pder{}{\theta_1} , r_2
    \pder{}{x_2} + \pder{}{\theta_2} \right\} .
\end{align*}
The expression for the projector $\P:T^*M \rightarrow C$ in local
coordinates is given by the following matrix
\begin{align*}
  \P = \left(
    \begin{array}{cccc}
      \frac{r_1^2}{1+r_1^2} & 0 & \frac{r_1}{1+r_1^2} & 0 \\
      0 & \frac{r_2^2}{1+r_2^2} & 0 & \frac{r_2}{1+r_2^2}\\
      \frac{r_1}{1+r_1^2} & 0 & \frac{1}{1+r_1^2} & 0\\
      0 & \frac{r_2}{1+r_2^2} & 0 & \frac{1}{1+r_2^2}
    \end{array}
  \right) .
\end{align*}
Let $\lambda \in T^*M_+$ be a plus-out point with $\ell (\lambda) = b$
and $\G (\lambda,d\ell)>0$. Since $H_- = \infty$, we only compute the
outcome of a reflective step. Following
equation~\eqref{eq:solutions-constrained-reflective}, the points in
the constrained characteristic passing through $\lambda$ with the same
$H_+$-energy level are $\lambda$ and $\lambda + c_{+,2} \P (d \ell)$,
with
\begin{align}\label{eq:the-constant}
  c_{+,2} = - \frac{2 \, \G(\lambda,\P(d\ell))}{\G(\P(d\ell),
    \P(d\ell))} \, .
\end{align}
According to Proposition~\ref{prop:duality}, the point $\lambda +
c_{+,2} \P (d \ell)$ is $+$-decisive.

Consider now an inelastic impact, i.e., the case when the length
$\ell$ of the rod becomes fixed after the impact. Since there are no
additional instantaneous constraints imposed on the system at the
impact state, we compute the decisive points with regards to the
boundary of the constraint manifold,
\begin{align*}
  C^\partial & = C \cap \alpha_{H_+}^{-1}(T\partial M) =
  \{(x_1,x_2,\theta_1,\theta_2,p_{x_1},p_{x_2},p_{\theta_1},p_{\theta_2})
  \in T^*M \; |  \; p_{x_1} = p_{x_2} ,\\
  & p_{x_1} = r_1 p_{\theta_1}, \, p_{x_2} = r_2 p_{\theta_2}, \,
  \ell(x_1,x_2,\theta_1,\theta_2) = a \; \text{or} \;
  \ell(x_1,x_2,\theta_1,\theta_2) = b\} \, .
\end{align*}
As before, we only compute the outcome of a reflective step. Following
Proposition~\ref{prop:decisive-dis-inelastic}, we deduce that the
unique decisive point is $\P^\partial(\lambda)$.  After the inelastic
impact, the length of the rod is fixed forever after, the velocities
of the two wheels of the system are reset according to
$\P^\partial(\lambda)$ and evolve according to
$X^{\partial}_{(H,C,N)}$.

\section*{Acknowledgments}
The work of the first author was partially supported by NSF grant
CMS-0100162.

\end{document}